%% file: Master.tex
\documentclass[12pt,reqno]{article}

\usepackage{cite,epic,eepic,euscript,verbatim,amsmath,amssymb,amsthm,amsfonts,afterpage,float,bm,stmaryrd,paralist}

\usepackage{cancel}
\usepackage{color}
\usepackage{fancyhdr}
\usepackage{graphics}              
\usepackage{hyperref}              
\usepackage{graphicx,epsf}
\usepackage{makeidx}
\usepackage{epsfig}
\usepackage{tikz}
\usepackage{subcaption}
\usetikzlibrary {decorations.markings,fpu} 

\newtheorem{theorem}{Theorem}[section]

\numberwithin{equation}{section}
\numberwithin{figure}{section}

\def\XXint#1#2#3{{\setbox0=\hbox{$#1{#2#3}{\int}$}
\vcenter{\hbox{$#2#3$}}\kern-.51\wd0}}

\newcommand{\reff}[1]{{\rm (\ref{#1})}}

\newcommand{\R}{\mathbb{R}}            
\newcommand{\ve}{\varepsilon}          

\graphicspath{{Figures/}}

\begin{document}

\title
{
Cell Polarity and Movement with Reaction-Diffusion and Moving Boundary: 
Rigorous Modeling and Robust Simulations
}

\author{
Shuang Liu
\thanks{Department of Mathematics, University of California, San Diego,
9500 Gilman Drive, La Jolla, California 92093-0112, United States.
Email: shl083@ucsd.edu
}
\and
Li-Tien Cheng
\thanks{
Department of Mathematics, University of California, San Diego,
9500 Gilman Drive, La Jolla, California 92093-0112, United States. 
Email: l3cheng@ucsd.edu
}
\and
Bo Li
\thanks{
Department of Mathematics and Quantitative Biology Ph.D. Program, 
University of California, San Diego,
9500 Gilman Drive, La Jolla, California 92093-0112, United States.
Email: bli@ucsd.edu
}
}

\date{\today}

\maketitle

\begin{abstract}

Cell polarity and movement are fundamental to many biological functions. 
Experimental and theoretically studies have indicated that interactions
of certain proteins lead to the cell polarization which plays a key role in 
controlling the cell movement. 
We study the cell polarity and movement based on a class of biophysical models
that consist of reaction-diffusion equations for different proteins and the 
dynamics of moving cell boundary. 
Such a moving boundary is often simulated by a phase-filed model. 
We first apply the matched asymptotic analysis to give a rigorous derivation
of the sharp-interface model of the cell boundary from a phase-field model.
We then develop a robust numerical approach that combines the level-set method
to track the sharp boundary of a moving cell and accurate discretization techniques
for solving the reaction-diffusion equations on the moving cell region.  
Our extensive numerical simulations predict the cell polarization under 
various kinds of stimulus, and capture both the linear and circular 
trajectories of a moving cell for a long period of time. 
In particular, we have identified some key parameters controlling 
different cell trajectories that are less accurately predicted by reduced models. 
Our work has linked different models and also developed tools that can be adapted 
for the challenging three-dimensional simulations.

\bigskip

\noindent
{\bf Key words and phrases}:
cell polarity, 
cell movement, 
reaction-diffusion equations, 
interface dynamics, 
matched asymptotic analysis, 
the level-set method.   
\end{abstract}

{\allowdisplaybreaks

\input{Introduction}

\input{Analysis}
\input{NumericalMethods}

\input{SimulationResults}

\input{Conclusions}


\bigskip

\noindent{\bf Acknowledgments.}
This work was supported in part by an AMS Simons Travel Grant (SL), 
the US National Science Foundation through grant DMS-1913144 (LTC $\&$ BL),
and the US National Institutes of Health through grant R01GM132106 (BL). 
The authors thank Dr.~Zirui Zhang and Professor Yanxiang Zhao for helpful discussions.

}

\bibliography{Cell}
\bibliographystyle{plain}

\end{document}

%% file: Introduction.tex
\section{Introduction}
\label{s:introduction}

Cell motility is fundamental to many biological functions such as immune response, morphogenesis, 
cancer metastasis, and wound healing, and yet it is extremely complicated \cite{Bray_CellMovement_2001}.
The movement of a eukaryotic cell crawling on a surface is a complex process, involving 
protrusion, retraction, and adhesion, exhibiting complex geometrical shapes and motion trajectories. 
Cell polarity, the spatially inhomogeneous distribution of different biomolecules 
such as proteins Rho GTPases inside a cell, resulting from biochemical interactions of biomolecules 
inside the cell, plays a crucial role in the cell movement 
\cite{Theriot_Nature2008,RappelKeshet_Rev2017,Mogilner_Science2012}. 
As cell polarity and movement involve multiple spatio-temporal scales and many-body interactions,
understanding such complex processes is challenging.

Recent years have seen the theoretical and computational development 
in studying cell polarity and motility \cite{Mogilner_Science2012,Mogilner_JMB2009,Keshet_Rev2020}.
Among different approaches, continuum models with reaction-diffusion equations and moving boundaries 
provide efficient simulation tools to understand the key mechanisms in 
cell polarity and movement 
\cite{rubinstein2005multiscale, 
vanderlei2011computational, 
camley2017crawling, 
mori2008wave,
Camley_PRL13,
Camley_PNAS14,
ShaoLevineRappel_PNAS12,
ShaoRappelLevine_PRL10,
Aranson_PLoSONE2013,
Keshet_SIAP2011}. 
An advantage of such modelling is that the motion of cell boundary, which is directly connected to the 
reaction and diffusion of a different biomolecules inside the cell, can be simulated
and analyzed to link the molecular basis for the 
cell polarity to the macroscopic cell movement, and to identify the key
parameters that control the cell polarization and movement.

In this work, we study the cell polarity and movement 
with a class of models that have been proposed in Mori et al.\
\cite{mori2008wave, Keshet_SIAP2011}, 
Shao et al.\ \cite{ShaoRappelLevine_PRL10},
and Camley et al.\ \cite{camley2017crawling,Camley_PRL13}. 
The basic components of such a model include a moving cell whose boundary motion is determined 
completely by its normal velocity, i.e., the normal component of the velocity. 
In addition to the geometrical effect, such normal velocity is controlled
by the amount of a membrane-bound Rho GTPases, e.g., Rac. 
Mori et al.\ \cite{mori2008wave} proposed the wave-pinning mechanism 
for the cell polarization: the reaction and diffusion of 
a membrane-bound active Rho GTPase and an inactive cytosolic form of Rho GTPase 
with a bistable kinetics lead to the formation of an 
interface inside a cell that separates a high from a low concentration of 
the active Rho GTPase proteins, and the propagation of such interface drives the cell polarization
that reaches a steady state eventually, as the wave is pinned down. 
Camley et al.\ \cite{camley2017crawling} reduce 
the two-species (active form and inactive form of Rho GTPase proteins) model proposed in
\cite{mori2008wave} to a single-species model, and also numerically simulated
the cell polarization and movement using a phase-field method, but carried out 
analysis with a sharp-interface description of the cell boundary motion. 
Simulations by Camley et al.\ \cite{camley2017crawling}
predict the linear and circular trajectories as a result of the wave-pinning dynamics.

To be specific, let us consider a moving cell confined spatially in a bounded region 
$\Omega \subset \R^d$ $(d = 2 \mbox{ or } 3)$.
Let us denote the cell boundary by $\Gamma(t)$ at time $t$ and assume it moves
with the normal velocity $V = V(x,t)$ for each point $x \in \Gamma(t)$. 
The cell boundary $\Gamma(t)$ separates the region $\Omega$ into the cell region, denoted
$\Omega^+(t)$, from the outer region, denoted $\Omega^-(t)$.  
As in \cite{camley2017crawling, mori2008wave},
let us consider the types of proteins inside the cell, a membrane-bound active and fast-diffusive
inactive Rho GTPase proteins, and denote their concentrations by $u = u(x,t)$ and $v = v(x,t)$, respectively. Extended from the one-dimensional model \cite{mori2008wave} (cf.\ also \cite{ShaoRappelLevine_PRL10}), 
our model of an underlying moving cell is
governed by the following system of equations and boundary conditions: 
\begin{align}
\label{NormalVelocity}
&
\tau V  = \alpha u  - \beta  -  \gamma H \quad \quad \mbox{for } x \in \Gamma(t) \mbox{ and } t > 0, 
\\ 
\label{RD4u}
&
\partial_t u = D_u \Delta u + f(u, v) \quad \mbox{for } x \in \Omega^+(t) \mbox{ and }  t > 0, 
\\
\label{RD4v}
&
\partial_t v = D_v \Delta v -  f(u, v) \quad  \mbox{for } x \in \Omega^+(t) \mbox{ and } t > 0, 
\\
\label{BC4uv}
&
\partial_n u = \partial_n v = 0 \ \qquad \qquad  \mbox{for } x \in \Gamma(t) \mbox{ and } t > 0, 
\end{align}
where 
\begin{equation}
\label{fForm}
f(u, v) = -k u ( u - 0.5 c) (u - C v). 
\end{equation}
In \reff{NormalVelocity}, $\tau$ is the friction coefficient, $\alpha$ and $\beta$ are 
the coefficients of F-actin extension and myosin retraction, respectively, 
$\gamma$ is the surface tension constant, and $H$ is the 
mean curvature of the cell boundary $\Gamma(t)$. In \reff{RD4u} and \reff{RD4v}, 
$D_u$ and $D_v$ are the diffusion constants for $u$ and $v$, respectively. 
We shall consider the regime that $D_v \gg D_u.$ In the reaction term $f(u, v)$ defined in \reff{fForm}, 
$k$ is the reaction rate relative to an average cell motility, 
$c$ is the constant averaged concentration $u$
in the cell front, and $C$ is a conversion parameter. 
Estimated values of these parameters are given in Table~1 in section~\ref{s:Results}. 
Note that the total mass
\begin{equation}
\label{M}
M = \int_{\Omega^+(t)} \left[ u(x,t) + v(x, t) \right] \, dx 
\end{equation}
is a constant with respect to time $t$.  

Assuming an infinite diffusion constant $D_v$, 
Camley et al.\ \cite{camley2017crawling}
have proposed and studied the following single-species model, reduced from the 
two-species model \reff{NormalVelocity}--\reff{BC4uv}: 
\begin{align}
\label{PFphi}
& \tau V = \alpha u - \beta - \gamma H 
\quad \, \qquad \mbox{for } x \in \Gamma(t) \mbox{ and }  t > 0, 
\\
\label{PFrho}
& \partial_t  u = D_u \Delta u + f(u, \bar{v})
\qquad \mbox{for } x \in \Omega^+(t) \mbox{ and } t > 0, 
\\
\label{PFphirhoBC}
& \partial_n u = 0 \,  \qquad  \qquad \qquad \qquad 
\mbox{for } x \in \Gamma(t)  \mbox{ and } t > 0,  
\end{align}
where $\bar{v}$, depending only on $u$,  is determined by the mass conservation and the 
fact that $v$ should be a constant with the large $D_v$ assumption, and is given by 
\[
\bar{v} = \frac{1}{\mbox{Area}\,(\Omega^+(t))} \left( M - \int_{\Omega^+(t)} u(x, t)\, dx\right). 
\]

To efficiently track the moving cell boundary in computer simulations, Shao et al.\ 
\cite{ShaoLevineRappel_PNAS12, ShaoRappelLevine_PRL10} and  
Camley et al.\ \cite{camley2017crawling, Camley_PRL13, Camley_PNAS14}
have used the phase-field model
\cite{
AndersonMcFaddenWheeler98, 
CollinsLevine85, 
Langer86, 
LRV_PRE09}.    
In such a model, 
the moving cell boundary is described by 
a continuous function, often called a phase field,
 that takes the value $1$ in the cell region and $0$ otherwise,
and smoothly changes its values from $0$ to $1$ in a thin transition layer, 
representing a diffuse cell boundary. 
Let us denote by $\phi_\ve = \phi_\ve(x,t)$ $(x \in \Omega, t \ge 0)$ such a phase-field function,
where $\ve \in (0, 1)$ is a small parameter and $t$ represents time.
Let us also denote by $u_\ve  = u_\ve (x,t)$ and $v_\ve = v_\ve (x, t)$ the
concentrations of the two different proteins, respectively, as described in \reff{RD4u} and \reff{RD4v}. 
Note that these functions are now defined on the entire region $\Omega.$ 
The phase-filed model that corresponds to the system of equations \reff{NormalVelocity}--\reff{BC4uv}, which 
shall be called a sharp-interface model, is then given by 
\begin{align}
\label{phiveIntro}
&
\tau \partial_t \phi_\ve  = (\alpha u_\ve - \beta)  | \nabla \phi_\ve | +
\gamma \left[ \Delta \phi_\ve  - \frac{1}{\ve^2} W'(\phi_\ve ) \right]
\quad \mbox{in } \Omega \times (0, \infty),
\\
\label{uveIntro}
&
\partial_t (\phi_\ve  u_\ve ) = \nabla \cdot D_u (\phi_\ve \nabla u_\ve ) + f(u_\ve, v_\ve )
\quad \mbox{in } \Omega \times (0, \infty),
\\
\label{vveIntro}
&
\partial_t (\phi_\ve  v_\ve ) = \nabla \cdot D_v (\phi_\ve  \nabla v_\ve ) - f(u_\ve, v_\ve )
\quad \mbox{in } \Omega \times (0, \infty),
\\
\label{BCIntro}
& \phi_\ve = u_\ve = v_\ve = 0 
\qquad  \mbox{in } \partial \Omega \times [0, \infty), 
\end{align}
where all the parameters and the function $f$ are the same as above,   
and $W = W(u)$ is a double-well potential given specifically by
\begin{equation}
\label{W}
W(u) = 18 u^2 (1 - u)^2 \qquad \forall u \in \R.
\end{equation}

In this work, we study the reaction-diffusion moving boundary model to understand the 
mechanisms of cell polarization and movement and the cooperation of these two processes. Our goal
is two fold. One is to understand the differences between some of the existing models
and make connections of such models. The other is to develop robust computational tools for
long-time accurate and efficient simulations of cell movement.  
Specifically: 
\begin{compactenum}
\item[(1)]
We derive rigorously the sharp-interface reaction-diffusion moving boundary model 
\reff{NormalVelocity}--\reff{BC4uv} from the phase-field model \reff{phiveIntro}--\reff{BCIntro}. 
While there are different advantages of different 
models, our analysis examines the consistency of these models. 
\item[(2)]
We develop a robust computational program that combines the level-set method and high-accurate
discretization method for solving reaction-diffusion equations on a moving cell region and 
for tracking the moving cell boundary. We test our numerical methods. 
\item[(3)]
We apply our numerical methods and algorithms to conduct a serious of computer 
simulations for the cell polarization and movement. We try to answer several questions: 
How does a cell respond to an external stimulus to polarize itself and then to move around? 
How does a cell keep different kinds of trajectories, such as a linear or a circular
trajectory, for a very long time?   
Our computational analysis predicts several important parameters such a finite diffusion 
constant (instead of taking it to be infinite in a reduced model), the surface tension
constant, and threshold concentration of an active Rho GTPase protein that partially controls the cell
movement. 
\end{compactenum}
Our computational tools prepare us well for future, large-scale three-dimensional simulations 
of the cell movement, which has been lacking currently in general. 

The paper is organized as follows. In section~\ref{s:SharpInterface}, 
we use the method of matched asymptotic analysis to derive 
the sharp-interface limit, the system \reff{NormalVelocity}--\reff{BC4uv}, 
of the phase-field reaction-diffusion moving boundary model
\reff{phiveIntro}--\reff{BCIntro}.
In section~\ref{s:NumericalMethods}, we describe a robust and accurate numerical method 
that combines a high-order finite difference discretization technique 
and a level-set method for the simulation of a moving cell. 
In section~\ref{s:Results}, we show our numerical simulations and analyze  our results
with various settings. Finally, in section~\ref{s:Conclusions}, 
we draw our conclusions and discuss several issues for further studies.

%% file: Analysis.tex
\section{From Phase-Field to Sharp-Interface Model}
\label{s:SharpInterface}

In this section, we carry out the matched asymptotic analysis 
\cite{Fife1988, Pego89, RubinsteinSternbergKeller89,LiLiu_SIAP2015,DaiPromislow13}
to derive the sharp-interface model \reff{NormalVelocity}--\reff{BC4uv} 
from the phase-field model \reff{phiveIntro}--\reff{BCIntro}. Specifically, we show  
that as $\ve \to 0$ the solution $\phi_\ve$ 
converges to the characteristic function of the cell region $\Omega^+(t)$, the normal velocity
of the cell boundary $\Gamma(t)= \partial \Omega^+(t)$  is given by \reff{NormalVelocity}, 
and the solutions $u_\ve$ and $\ve_\ve$ converge to the solutions to \reff{RD4u}--\reff{BC4uv}.  

We shall analyze the following more general phase-field model
in the setting of three-dimensional space: 
\begin{align}
\label{phive}
&
\partial_t \phi_\ve  = h(u_\ve, v_\ve)| \nabla \phi_\ve | + \gamma \left[ \Delta \phi_\ve  
- \frac{1}{\ve^2} W'(\phi_\ve ) \right] \quad \mbox{in } \Omega \times (0, \infty),
\\
\label{uve}
&
\partial_t (\phi_\ve u_\ve ) = \nabla \cdot D_1 (\phi_\ve \nabla u_\ve) + f(u_\ve, v_\ve)
\quad \mbox{in } \Omega \times (0, \infty),
\\
\label{vve}
&
\partial_t (\phi_\ve v_\ve ) = \nabla \cdot D_2 (\phi_\ve \nabla v_\ve) + g(u_\ve, v_\ve)
\quad \mbox{in } \Omega \times (0, \infty),
\\
\label{BC}
& \phi_\ve = u_\ve = v_\ve = 0 \qquad  \mbox{in } 
\partial \Omega \times [0, \infty).
\end{align}
Here, $\Omega \subset \R^3$ is a smooth and bounded domain, 
$\ve \in (0, 1)$ is a small parameter,
$\gamma > 0$, $D_1 > 0$, and $D_2 > 0$ are all constants, 
and $f$, $g$, and $h$ are all smooth and bounded two-variable functions. 
The double-well function $W $ is defined in \reff{W}. 
Note that the analysis for the single-species system \reff{PFphi}--\reff{PFphirhoBC} or for 
a two-dimensional setting is similar.

\medskip

\noindent
{\bf Initial formation of a diffuse cell boundary.}
We assume the following expansions:
\begin{align*}
&\phi_\ve (x,t) = \phi_0(x, \tau) + \ve \phi_1 (x, \tau) + \ve^2 \phi_2(x, \tau) + \cdots,
\\
&u_\ve (x,t) = u_0(x, \tau) + \ve u_1 (x, \tau) + \ve^2 u_2(x, \tau) + \cdots,
\\
&v_\ve (x,t) = v_0(x, \tau) + \ve v_1 (x, \tau) + \ve^2 v_2(x, \tau) + \cdots,
\end{align*}
where $\tau = \tau(t, \ve) $ is a time variable that can be 
different from the regular time variable $t$,
and all the functions $\phi_i = \phi_i(x,\tau)$, $u_i = u_i (x, \tau)$, and
$v_i = v_i(x,\tau) $  $(i = 0, 1, \dots) $ are smooth and bounded in $\Omega,$ satisfying the 
boundary conditions $\phi_i = u_i = v_i = 0$ on $\partial \Omega;$ cf.\ \reff{BC}.

Considering a fast time scale $\tau = t/\ve^2$, we have $\partial_t = \ve^{-2} \partial_\tau$.
Plugging the above expressions of $\phi_\ve$, $u_\ve$, and $v_\ve$ into \reff{phive}, 
using Taylor's expansion, and comparing terms of the leading orders $O(\ve^{-2})$ 
and $O(\ve^{-1})$, respectively, we obtain that
\[
O(\ve^{-2}): \quad \partial_\tau \phi_0 = - \gamma W'(\phi_0) \qquad \mbox{and} \qquad
O(\ve^{-1}): \quad \partial_\tau \phi_1  = - \gamma W''(\phi_0) \phi_1.
\]
Since $W'(s) = 0$ if and only if $s = 0$, $1/2$, or $1$, with $0$ and $1$ being local minima
of $W$ and $1/2$ being a local maximum of $W$,
 given any point $x\in \Omega$ and any initial data $\phi_0(x,0) \ne 1/2$
$\phi_0(x,\tau) \to 0 $ or $1$ exponentially as $\tau \to \infty.$
Once $\phi_0$ falls into $(-\infty, (3-\sqrt{3})/6) \cup ( (3 + \sqrt{3})/6, \infty),$
the convex region of $W$, then $\phi_1(x, \tau) \to 0$ exponentially as $\tau \to \infty.$
If we consider the next fast time scale with $\tau = t/\ve$, then
we have $\partial_t = \ve^{-1} \partial_\tau.$ Similar calculations lead to
the leading-order equations
\[
O(\ve^{-2}): \quad W'(\phi_0) = 0 \qquad \mbox{and} \qquad O(\ve^{-1}): \quad
\partial_\tau \phi_0 = - \gamma  W''(\phi_0) \phi_1.
\]
Again, we see that $\phi_0 = 0$, $1/2$, or $1$, and since $W''(\phi_0) \ne 0$, we have $\phi_1 = 0.$
Results are the same if we consider the regular time scale $\tau = t.$

We can therefore 
assume that the region $\Omega$ is divided by the phase-field function $\phi_\ve$ into
an outer region $O_\ve(t) := \Omega_\ve^-(t) \cup \Omega_\ve^+(t)$, where
\[
 \Omega_\ve^- (t) = \{ x\in \Omega: u_\ve(x,t) = O(\ve^2)  \}
\quad \mbox{and} \quad 
 \Omega_\ve^+(t) = \{ x\in \Omega: u_\ve(x,t) = 1 + O(\ve^2) \},
\]
and an inner region $I_\ve(t) := \Omega \setminus O_\ve(t)$, 
 where $\phi_\ve$ changes from $0$ to $1$,  
representing the diffuse cell boundary.
The region $\Omega_\ve^+(t)$ is the cell region at $t$. 
Note that, by the imposed boundary conditions $u_\ve = 0 $ on $\partial \Omega$,
the boundary $\partial \Omega$ is included in the closure of $\Omega_\ve^-(t).$
We further assume that the inner region $I_\ve(t)$
is an $O(\ve)$-neighborhood of a closed and smooth surface $\Gamma(t)$, independent of $\ve$,
that is the limit of $\{ x \in \Omega: u_\ve(x, t) = 1/2 \}$ as $\ve \to 0.$
Moreover, the interior and exterior of $\Gamma(t)$, denoted $\Omega^+(t)$  
and $\Omega^-(t)$, are the limit as $\ve \to 0$ of $\Omega^+_\ve(t)$ and $\Omega_\ve^-(t)$, respectively, 
with $\Omega^+(t)$ being the cell region.

\medskip

\noindent
{\bf Outer expansions.}
We assume the following expansions in the outer region $\Omega_\ve(t)$: 
\begin{align*}
&\phi_\ve (x,t) = \phi_0(x, t) + \ve \phi_1 (x, t) + \ve^2 \phi_2(x, t) + \cdots,
\\
&u_\ve (x,t) = u_0(x, t) + \ve u_1 (x, t) + \ve^2 u_2(x, t) + \cdots,
\\
&v_\ve (x,t) = v_0(x, t) + \ve v_1 (x, t) + \ve^2 v_2(x, t) + \cdots,
\end{align*}
where the functions $\phi_i(x, t)$, $u_i(x,t),$ and $v_i(x,t)$ $(i = 0, 1, \dots)$ are
smooth and bounded, and are independent of $\ve$. They also satisfy the boundary conditions
$\phi_i = u_i = v_i = 0$ on $\partial \Omega.$  Note that these functions are
different from those in the expansions with a different time scale $\tau.$
Since $\phi_0=O(\ve^2)$ in $\Omega_\ve^{-}(t)$, there will be
no equations for $u_\ve$ and $v_\ve$ at leading order $O(1)$, we shall assume that
$ u_\ve = 0$ and $v_\ve = 0$ in $\Omega_\ve^-(t).$
If we plug the above expansion of $\phi_\ve$, $u_\ve$, and $v_\ve$
into \reff{phive}, \reff{uve}, and \reff{vve},
we obtain by a series of calculations that, up to the leading order $O(1),$
\begin{align}
\label{u0outer}
& \partial_t u_0 = D_1 \Delta u_0 + f(u_0, v_0)
\quad \mbox{in } \Omega^+(t) \times (0, \infty),
\\
\label{v0outer}
& \partial_t v_0  =  D_2 \Delta v_0 + g(u_0, v_0)
\quad \mbox{in } \Omega^+(t)+ \times (0, \infty).
\end{align}

\medskip

\noindent
{\bf Local coordinates for the inner region.}
Let $x \in I_\ve(t)$ and denote by $s(x,t)$ the signed distance from $x$ to $\Gamma(t)$,
with $s (x,t)> 0$ if $x$ is inside $\Gamma(t)$ and $s (x,t) < 0$ otherwise.
Note that $s(x,t) = O(\ve)$ and $|\nabla s (x,t) | = 1.$
Now, let $y = P(x,t) \in \Gamma(t)$ be the projection of $x$ onto $\Gamma(t)$, defined by
$| x - P(x,t) | = |s(x,t) | $.  (We use $| \, \cdot \, | $ to denote both the
absolute value of a number and the Euclidean norm of a vector.)
For $0 < \ve \ll 1$, the projection $y = P(x,t) \in \Gamma(t)$ is unique, and
the vector $x - P(x,t) $ is normal to the surface $\Gamma(t)$ at $y = P(x,t)$.
Let $z = s(x,t)/\ve$.
Let $n = n(y,t) = \nabla s(y,t)$ be the unit normal at $y \in \Gamma(t)$ pointing from
the exterior to the interior of $\Gamma(t)$.
Then we have a unique expression of $x \in I_\ve(t)$ as
\begin{equation}
\label{xyz}
x = y + \ve z  n.
\end{equation}
We call $(y,z)$ the local coordinate of $x \in I_\ve(t)$ with respect to the surface $\Gamma(t).$
We have for $0 < \ve \ll 1$ that
\cite{Fife1988,Pego89, RubinsteinSternbergKeller89,LiLiu_SIAP2015,DaiPromislow13}
\begin{align}
\label{dxz}
& \nabla_x z = \ve^{-1} n(y,t) + O(1),
\\
\label{ddxz}
& \Delta_x z = 2 \ve^{-1} H(y,t) + O(1),
\\
\label{dtz}
& \partial_t z = - \ve^{-1} V(y,t),
\\
\label{dyn0}
& \nabla_x y_j(x, t) \cdot n(y,t) = 0  \quad (j = 1, 2, 3),
\end{align}
where
$H(y,t)$ is the mean curvature of the surface $\Gamma(t)$ at the point $y = P(x,t)$,
$V(y,t)$ is the  normal velocity of the point $y = P(x,t) \in \Gamma(t)$ defined by
\begin{equation}
\label{vyt}
V(y,t) = \partial_t y \cdot n(y,t) = \partial_t P(x,t) \cdot n(y,t),
\end{equation}
and $y_j$ $(j=1,2,3)$ are the components of $y = y(x,t)$.
Let  $f = f(x,t)$ and $\tilde{f}  = \tilde{f}(z,y,t)$ be smooth functions such that 
$f(x,t) = \tilde{f}(z,y,t)$  with $x \in I_\ve(t)$ and $(y, z)$ related by \reff{xyz}. Then,
by \reff{dxz}--\reff{vyt} and the chain rule,
we obtain for $0 < \ve \ll 1$ that
\cite{Fife1988,Pego89, RubinsteinSternbergKeller89,LiLiu_SIAP2015,DaiPromislow13}
\begin{align}
\label{x1}
& \nabla_x f (x,t) =  \ve^{-1} n \, \partial_z  \tilde{f}(y, z, t) + O(1 ),\\
\label{x2}
&\Delta_x f (x, t) =  \left( \ve^{-1} 2 H (y,t) \partial_z +\ve^{-2}\partial^2_{zz} \right)
\tilde{f}(y, z, t) +O(1), \\
\label{x3}
& \partial_t f (x, t) = - \ve^{-1} V(y,t) \, \partial_z \tilde{f}(y, z, t) + O(1).
\end{align}

\medskip

\noindent
{\bf Inner expansions.}
We now assume that following expansions 
in the inner region $I_\ve(t)$: 
\begin{align*}
& \phi_\ve(x,t) = \tilde{\phi}_0(y,z, t) + \ve\tilde{\phi}_1(y,z, t)
+ \ve^2\tilde{\phi}_2(y,z, t)+\cdots,
\\
& u_\ve (x,t) = \tilde{u}_0 (y,z, t)+ \ve \tilde{u}_1(y,z, t) +
\ve^2\tilde{u}_2 (y,z, t) + \cdots,
\\
& v_\ve (x,t) = \tilde{v}_0 (y,z, t)+ \ve \tilde{v}_1(y,z, t) +
\ve^2\tilde{v}_2 (y,z, t) + \cdots,
\end{align*}
where $x \in I_\ve(t)$ and $(y,z)$ are related by \reff{xyz},
and all $\tilde{\phi}_i = \tilde{\phi}_i(y,z, t), $  $\tilde{u}_i=  \tilde{u}_i(y, z, t)  $, and
$\tilde{v}_i = \tilde{v}_i(y, z, t)$  $(i = 0, 1, \dots)$ are smooth and bounded functions.
Let us substitute $\phi_\ve$, $u_\ve$, and $v_\ve$ in \reff{phive}--\reff{vve} 
with these expansions. By \reff{x1}--\reff{x3} and a series of calculations, we obtain
\begin{align*}
- \ve^{-1} V \partial_z \tilde{\phi}_0 &= \ve^{-1} h(u_0, v_0)
\left| \partial_z \tilde{\phi}_0 \right|
+ \ve^{-2} \gamma \left[ \partial_{zz} \tilde{\phi}_0 - W'(\tilde{\phi}_0) \right]
\\
& \qquad
+ \ve^{-1} \gamma \left[ 2 H \partial_z \tilde{\phi}_0 +
\partial_{zz} \tilde{\phi}_1 - W''(\tilde{\phi}_0) \tilde{\phi}_1 \right] + O(1),
\end{align*}
where $V = V(y,t)$ and $H = H(y,t)$ are the normal velocity and mean curvature, respectively,
at $y = P(x,t).$
Note that, unlike $\phi_\ve$ which varies from $0$ to $1$,
the concentration fields $u_\ve$ and $v_\ve$ should not vary largely, in the inner region.
In particular, we have
 $ \tilde{u}_0 (y, z, t) = \tilde{u}_0(y, 0, t)+O(\ve)$ and 
 $ \tilde{v}_0 (y, z, t) = \tilde{v}_0(y, 0, t)+O(\ve),$ 
as $ |x - y | = | x - P(x,t)| = O(\ve)$ for any $x \in I_\ve(t)$ with the local coordinate $(y,z)$.
Therefore, 
we obtain from the above equation that
\begin{align*}
- \ve^{-1} V \partial_z \tilde{\phi}_0
&= \ve^{-1}  h(\tilde{u}_0(y, 0, t), \tilde{v}_0(y, 0, t))
\left| \partial_z \tilde{\phi}_0 \right|
+ \ve^{-2} \gamma \left[ \partial_{zz} \tilde{\phi}_0 - W'(\tilde{\phi}_0) \right]
\\
& \qquad
+ \ve^{-1} \gamma \left[ 2 H \partial_z \tilde{\phi}_0 +
\partial_{zz} \tilde{\phi}_1  - W''(\tilde{\phi}_0) \tilde{\phi}_1 \right] + O(1).
\end{align*}
Now, equating the terms with the same order $O(\ve^{-2})$ and $O(\ve^{-1})$,
respectively, we get
\begin{align}
\label{tildephi0}
&O(\ve^{-2}): \qquad  0 = \partial_{zz} \tilde{\phi}_0 - W'(\tilde{\phi}_0),
\\
\nonumber
& O(\ve^{-1}): \qquad - V \partial_z \tilde{\phi}_0 =
 h(\tilde{u}_0(y, 0, t), \tilde{v}_0(y, 0, t))
\left| \partial_z \tilde{\phi}_0 \right|
\\
\label{phi0phi1}
&\qquad \qquad \qquad \qquad \qquad +  \gamma \left[ 2 H \partial_z
\tilde{\phi}_0 + \partial_{zz} \tilde{\phi}_1
- W''(\tilde{\phi}_0) \tilde{\phi}_1 \right].
\end{align}
Similarly, we can plug the inner expansions of $\phi_\ve$, $u_\ve$, and $v_\ve$
into \reff{uve} and \reff{vve} to get
in the leading order that
\begin{align}
\label{tildephi0u0}
& O\left(\ve^{-2}\right): \qquad 0 = D_1 \left( \partial_z \tilde{\phi}_0
\partial_z \tilde{u}_0 + \tilde{\phi}_0 \partial_{zz} \tilde{u}_0  \right),
\\
\label{tildephi0v0}
& O\left(\ve^{-2}\right): \qquad 0 = D_2 \left( \partial_z \tilde{\phi}_0
\partial_z \tilde{v}_0 + \tilde{\phi}_0 \partial_{zz} \tilde{v}_0  \right).
\end{align}

\medskip

\noindent
{\bf Inner-outer matching and the sharp-interface limit.}
Since in the outer region $\phi_\ve = O(\ve^2) $ in $\Omega_\ve^-(t)$
and $\phi_\ve = 1+O(\ve^2) $ in $\Omega_\ve^+(t)$, we have the following matching conditions for the leading-order terms of
the inner and outer solutions of the phase field $\phi_\ve$:
\begin{equation}
\label{matchu0}
\lim_{z \to -\infty} \tilde{\phi}_0 (y, z, t) = 0  \qquad
\mbox{and} \qquad
\lim_{z \to \infty} \tilde{\phi}_0 (y, z, t) = 1.
\end{equation}
These, together with \reff{tildephi0}, determine completely $\tilde{\phi}_0$ to be
\begin{equation*}
\tilde{\phi}_0(y,z, t) =
\frac12 + \frac{e^{3z}-e^{-3z}}{2(e^{3z}+e^{-3z})} \qquad \forall z \in \R.
\end{equation*}
In particular, $\tilde{\phi}_0$ does not depend on $y$ and $t$.
One can verify that $\partial_z \tilde{\phi}_0 > 0$ and that  
\begin{equation}
\label{one}
 \int_{-\infty}^{\infty} (\partial_z \tilde{\phi}_0)^2{\rm d}z
= 1. 
\end{equation}
By matching the inner expansion and outer expansion, we have
$\partial_z  \tilde{\phi}_0(\pm \infty) = \partial_{zz} \tilde{\phi}_0(\pm \infty) = 0$.
Thus, by integration by parts and \reff{tildephi0}, we have 
\begin{equation}
\label{nice}
\int_{-\infty}^{\infty} \partial_z \tilde{\phi}_0 \left[  \partial_{zz}\tilde{\phi}_1
- W''(\tilde{\phi}_0)\tilde{\phi}_1 \right] dz =  \int_{-\infty}^{\infty} \partial_z
\left[ \partial_{zz} \tilde{\phi}_0 - W'(\tilde{\phi}_0 ) \right] \tilde{\phi}_1 \, {\rm d}z = 0.
\end{equation}
Now, by multiplying both sides of \reff{phi0phi1}
by $\partial_z \tilde{\phi}_0$ and then integrating the resulting equation
over $z \in (-\infty, \infty)$, we have by \reff{one} and \reff{nice} that 
\begin{equation}
\label{VVV}
 V(y,t) = - h(u_0(y, 0, t), v_0(y, 0, t)) - 2 \gamma H (y, t) \qquad \forall y \in \Gamma(t).
\end{equation}

It follows from \reff{tildephi0u0} and \reff{tildephi0v0} that
$\partial_z ( \tilde{\phi}_0 \partial_z \tilde{u}_0) = 0$ and
$\partial_z ( \tilde{\phi}_0 \partial_z \tilde{v}_0) = 0$ for all $z \in \R.$ These
and the matching conditions \reff{matchu0} imply that
\begin{align*}
&0 = \lim_{z \to - \infty} \tilde{\phi}_0(y, z, t) \partial_z \tilde{u}_0(y, z, t) =
\lim_{z \to \infty} \tilde{\phi}_0(y, z, t) \partial_z \tilde{u}_0(y, z, t)
=  \lim_{z \to \infty} \partial_z \tilde{u}_0(y, z, t),
\\
&0 = \lim_{z \to - \infty} \tilde{\phi}_0(y, z, t ) \partial_z \tilde{v}_0(y, z, t) =
\lim_{z \to \infty} \tilde{\phi}_0(y, z, t ) \partial_z \tilde{v}_0(y, z, t)
= \lim_{z \to \infty} \partial_z \tilde{v}_0(y, z, t).
\end{align*}
Matching the inner expansions and outer expansions, we thus have
\begin{align}
\label{u0Neumann}
&\partial_n u_0(x, t) = \lim_{z\to \infty} \partial_z \tilde{u}_0(y, z, t) = 0
\qquad \forall x \in \Gamma(t),
\\
\label{v0Neumann}
& \partial_n v_0(x, t) = \lim_{z\to \infty} \partial_z \tilde{u}_0(y, z, t) = 0
\qquad \forall x \in \Gamma(t),
\end{align}
where $u_0$ and $v_0$ are the leading-order terms in the outer expansion of $u_\ve$ and
$v_\ve$, and satisfy \reff{u0outer} and \reff{v0outer}, respectively.

We summarize our analysis in the following:
\begin{theorem}
\label{t:SharpLimit}
Under the assumption that there exists a closed and smooth interface $\Gamma(t)$,
the outer and inner expansions above for the solutions
$\phi_\ve$, $u_\ve, $ and $v_\ve$ are valid, and 
the corresponding matching conditions are satisfied, the following hold true
in the limit $\ve \to 0$: 
\begin{compactenum}
\item[\rm{(1)}]
The  phase-field function $\phi_\ve$ converges to
$1$ in $\Omega^+(t)$ and $0$ in $\Omega^-(t)$, respectively;
\item[\rm{(2)}]
The concentrations $u_\ve$ and $v_\ve$ converge to the 
solution to the boundary-value problem of the reaction-diffusion equations
\reff{u0outer} and \reff{v0outer}, and \reff{u0Neumann} and \reff{v0Neumann}; 
\item[\rm{(3)}]
The normal velocity $V=V(x,t)$ 
of the sharp cell boundary $\Gamma(t)$ is given by \reff{VVV}. 
\end{compactenum}
\end{theorem}

%% file: NumericalMethods.tex
\section{Numerical Methods}
\label{s:NumericalMethods}

We describe our numerical methods for solving the system of equations
\reff{NormalVelocity}--\reff{BC4uv} in the following non-dimensionalized form
for the rescaled normal velocity $V = V(x, y, t)$ on the rescaled cell boundary 
$\Gamma(t)$ at time $t$, rescaled concentrations
$u = u(x, y, t)$ and $v = v(x, y, t)$ defined on the cell region $\Omega^+(t)$, respectively: 
\begin{align}
\label{eq:nonV}
&V =  u-u^*  - \chi H \qquad \ \ \, \quad \text{for } (x, y) \in \Gamma(t) \mbox{ and } t > 0, 
\\
\label{eq:nonEqu}
& 
\partial_{t}{u}=D_u\Delta u+f(u,v) \qquad \text{for }  (x, y) \in \Omega^+(t) \mbox{ and } t > 0, 
\\
\label{eq:nonEqv}
&
\partial_{t}{v}=D_v\Delta v-f(u,v) \qquad \text{for } (x, y) \in \Omega^+(t) \mbox{ and } t > 0, 
\\
\label{eq:nonBC}
&
\partial_n u=\partial_n v=0 \qquad \quad \, \qquad \text{for } (x, y) \in \Gamma(t) \mbox{ and } t > 0,  
\end{align}
where $H$ is the rescaled curvature, all  $u^\ast$, $\chi$, 
$D_u$, $D_v$, $K$, and $C$ are positive constants, and  
\[
f(u,v)=-Ku(u-0.5)(u-Cv).  
\]
Details of the non-dimensionalization are given in section~4.1 below. 
We solve these equations in two-dimensional space.   
Our numerical methods for the one-species system, which is the sharp-interface limit of the system 
\reff{PFphi}--\reff{PFphirhoBC}, are similar. 

We set our computational box to be $\Omega = (-L, L)^2$ for some $L>0$ and cover it with a uniform
finite-difference grid with stepsize $h$ in each dimension. We discretize a time 
interval $[0, T]$
with $T > 0$ the final time of interest by $t_m = m \Delta t$ 
$(m = 0, 1, \dots)$ with time step $\Delta t > 0$. 
We denote by $\Gamma_m$, $\Omega^+_m$, and $\Omega^-_m$ the approximation
of $\Gamma(t_m)$, $\Omega^+(t_m)$, and $\Omega^-(t_m)$, respectively, 
where $\Omega^-(t)$ is the rescaled outer region. 
For a function $w = w(x,y, t)$ with  $ (x, y) \in \Omega$ and $t \ge 0$, 
we denote by $w^m = w^m(x,y)$ the approximation of $w(x,y, t_m)$,    
and by $w^m_{i,j}$ the approximation of $w^m(x_i,y_j)$ for a grid point $(x_i,y_j)$ 
or the center $(x_i, y_j)$ of a grid cell. 
Note that we approximate the level-set function $\phi$ at grid points while we 
approximate the concentrations $u$ and $v$ at the centers of grid cells; see below. 

\medskip

\noindent
{\bf The Level-set method for the moving cell boundary.}
We capture the cell boundary $\Gamma(t)$ at time $t$ by using the level-set method 
\cite{OS88,OF02}, with level-set function $\phi = \phi(x,y, t),$
i.e., $\Gamma(t) = \{ (x, y) \in \Omega: \phi(x,y, t) = 0 \}.$    
The level-set function is determined by the evolution equation
$\partial_t \phi + V |\nabla \phi | = 0,$
where $V = V(x,y, t)$ is given in \reff{eq:nonV} that needs to be 
 extended from $\Gamma(t)$ to the entire  computational domain $\Omega.$ 
The first part of our
normal velocity is $u (x, y, t) - u^\ast$. We keep the value of $u = u(x,y, t)$ in 
$\Omega^+(t)$ and additionally extend it from $\Gamma(t)$ to $\Omega^-(t)$ numerically 
in each step of time iteration.  Note, for convenience, we continue to denote the result 
by $ u = u(x,y, t)$, now for $(x, y) \in \Omega$. (Details of such extension are given below.)
The curvature $H$  can be extended simply 
by using $H  = \nabla \cdot ( \nabla \phi / | \nabla \phi | )$ for all $(x, y)\in\Omega$. 
Note that implicitly we require that the level-set function to be close to the signed distance
to the interface $\Gamma(t)$ with $\phi < 0$ in $\Omega^+(t)$ (the cell region)
and $\phi > 0$ in $\Omega^-(t)$, 
at least near $\Gamma(t)$ Therefore, the level-set equation and boundary conditions become now 
\begin{align}
\label{eq:level_set_cell_motion1} 
& \partial_{t} \phi =  -(u-u^*)  | \nabla \phi| +\chi
\big(\nabla \cdot \frac{ \nabla \phi}{| \nabla \phi|}\big)  | \nabla \phi| 
& &  \mbox{for } (x, y) \in \Omega \mbox{ and }  t > 0,  & 
\\     
\label{eq:level_set_cell_motion2}
&  \partial_n \phi = 0 & & \mbox{for } (x, y) \in \partial \Omega \mbox{ and } t > 0. 
\end{align}
Following \cite{smereka2003semi}, we rewrite the curvature part of the normal velocity as 
\begin{align*}
& \nabla \cdot\left( \frac{ \nabla \phi}{| \nabla \phi|} \right)
 | \nabla \phi| =\Delta \phi -N(\phi),
\\
&     N(\phi)=\frac{\nabla \phi}{|\nabla \phi|}\cdot \nabla (|\nabla \phi|)
= \frac{\phi_x^2\phi_{xx}+2\phi_x\phi_y\phi_{xy}+\phi_y^2\phi_{yy}}{\phi_x^2+\phi_y^2}. 
\end{align*}
With a given initial level-set function $\phi(x,y, 0)$ for all $(x, y) \in \Omega$, we can solve 
\reff{eq:level_set_cell_motion1} and \reff{eq:level_set_cell_motion2} numerically 
with finite difference schemes in time and over the uniform grid.

Specifically, starting from $u^m_{i,j}$ at centers of grid cells, we first use 
polynomial interpolations or extrapolations to approximate $u$ at points where
the interface intersects grid lines, itself computed using linear interpolation 
on values of $\phi$ at grid points.  We then extend these values from the points on 
the interface to all the grid points 
in the outer region $\Omega_m^+$, constant in the normal direction, by the fast sweeping 
method\cite{Z05,tsai2003fast}.

We continue to denote the extended function by $u^m$. 
To then get $\phi^{m+1}$, we use the semi-implicit scheme 
\begin{equation}
\label{phiphi}
        \frac{\phi^{m+1}-\phi^{m}}{\Delta t}
=- (u^m - u^\ast) | \nabla \phi^m |  + \chi  \Delta \phi^{m+1} - \chi  N_\epsilon(\phi^m),  
\end{equation}
where $N_\epsilon (\phi)$ is the same as $N(\phi^m)$ except the denominator 
$\phi_x^2+\phi_y^2$ in $N(\phi)$ is replaced by $\phi_x^2+\phi_y^2+\epsilon$ 
in $N_\epsilon(\phi)$ for a small enough $\epsilon > 0$, 
to avoid singularities at $\nabla\phi = 0$ while keeping an accurate approximation away from them. 
We discretize $ | \nabla \phi^m |$ by fifth-order WENO\cite{JP00} within
Godunov's scheme\cite{OS91}, and discretize 
$\Delta \phi$ and $N_\epsilon(\phi)$ by second-order central differencing.

At the boundary of the computational domain $\Omega$, 
we use a second-order scheme to discretize the Neumann boundary conditions
\reff{eq:level_set_cell_motion2}.  
The coefficient matrix of the resulting system of linear equations for all 
$\phi^{m+1}_{i,j}$ is sparse and nonsymmetric, with the nonsymmetry due mainly to the 
chosen treatment of the boundary conditions.
We solve the linear system of equations using the 
biconjugate gradient stabilized method preconditioned with the incomplete LU decomposition. 
Finally, we reinitialize the level-set function $\phi^{m+1}$, performing a few iterations
of the algorithm of redistancing to signed distance function \cite{SSO94}, and 
continue to denote the result by $\phi^{m+1}.$

\medskip

\begin{figure}[H]
	\centering 
	\begin{tikzpicture} [line width=0.5, scale=1.2] 
		\clip (-1.1,-1.1) rectangle (5,4);
		\draw [blue, fill=blue!100, opacity=0.8] plot [smooth cycle] coordinates {(1.0,.1)(1.75,0.05)(2.5,.5)(2.9,1.5)(2.8,2.8)(2.0, 2.9)(1.0,2.6)(0.3, 1.3)(0.5,0.5)};
		\draw [red, fill=red!80, opacity=0.8] plot [smooth cycle] coordinates { (3.3, 2.55)(3.2, 1.8)(2.5, 0.7) (1.0,0.4) (0.6,1.8)(1.2, 3.15)(2.65, 3.35)};
		
		\draw [blue] plot [smooth cycle] coordinates {(1.0,.1)(1.75,0.05)(2.5,.5)(2.9,1.5)(2.8,2.8)(2.0, 2.9)(1.0,2.6)(0.3, 1.3)(0.5,0.5)};
		\node (A) at (3.3, 3.3) {$\Gamma_{m+1}$};
		\node (A) at (1.0,-0.2) {$\Gamma_{m}$};
		
		\node (A) at (1.85,0.28) {III};
		
			\node (A) at (1.7,2.0) {II};
			
				\node (A) at (1.6,3.1) {I};
		
		\filldraw [black] (0.7,2.5) circle (0.5pt);
		\filldraw [black] (0.9,2.7) circle (0.5pt);
		\filldraw [black] (1.1,2.9) circle (0.5pt);
		\filldraw [black] (0.9,2.9) circle (0.5pt);
		\filldraw [black] (1.1,3.1) circle (0.5pt);
			\filldraw [black] (1.3,2.9) circle (0.5pt);
		\filldraw [black] (1.3,3.1) circle (0.5pt);
		\filldraw [black] (1.3,3.3) circle (0.5pt);
		\filldraw [black] (1.5,3.1) circle (0.5pt);
		\filldraw [black] (1.5,3.3) circle (0.5pt);
		\filldraw [black] (1.7,3.3) circle (0.5pt);
			\filldraw [black] (1.7,3.1) circle (0.5pt);
				\filldraw [black] (1.9,3.1) circle (0.5pt);
					\filldraw [black] (1.9,3.3) circle (0.5pt);
						\filldraw [black] (2.1,3.1) circle (0.5pt);
		\filldraw [black] (2.1,3.3) circle (0.5pt);
		\filldraw [black] (2.3,3.1) circle (0.5pt);
		\filldraw [black] (2.3,3.3) circle (0.5pt);
		\filldraw [black] (2.5,3.3) circle (0.5pt);
		\filldraw [black] (2.5,3.1) circle (0.5pt);
		\filldraw [black] (2.7,3.1) circle (0.5pt);
		\filldraw [black] (2.7,3.3) circle (0.5pt);
		\filldraw [black] (2.9,3.1) circle (0.5pt);
		\filldraw [black] (2.9,3.3) circle (0.5pt);
			\filldraw [black] (2.9,2.9) circle (0.5pt);
				\filldraw [black] (3.1,3.1) circle (0.5pt);
					\filldraw [black] (3.1,2.9) circle (0.5pt);
						\filldraw [black] (3.1,2.7) circle (0.5pt);
							\filldraw [black] (3.1,2.5) circle (0.5pt);
								\filldraw [black] (3.1,2.3) circle (0.5pt);
									\filldraw [black] (3.1,2.1) circle (0.5pt);
										\filldraw [black] (3.1,1.9) circle (0.5pt);
											\filldraw [black] (3.3,1.9) circle (0.5pt);
												\filldraw [black] (3.3,2.1) circle (0.5pt);
													\filldraw [black] (3.3,2.3) circle (0.5pt);
														\filldraw [black] (3.3,2.5) circle (0.5pt);
															\filldraw [black] (3.3,2.7) circle (0.5pt);
																\filldraw [black] (3.1,1.7) circle (0.5pt);
																	\filldraw [black] (3.1,1.5) circle (0.5pt);
																	\filldraw [black] (3.1,1.3) circle (0.5pt);
		
		\draw[step=.2cm,gray, very thin] (-1,-1) grid (4,4);
	\end{tikzpicture}
	\caption{Illustration of two consecutive cell regions $\Omega_m^+$ and $\Omega_{m+1}^+$. 
The cell region $\Omega_m^+$ is the union of part II and part III, enclosed by the cell 
boundary $\Gamma_m$ (blue curve). 
The cell region $\Omega_{m+1}^+$ is the union of part I and part II, marked by red. 
Black solid dots mark those centers of grid cells that are in $\Omega_{m+1}^+$ but not in 
$\Omega_m^+.$}
\label{fig:cell_move}
\end{figure}
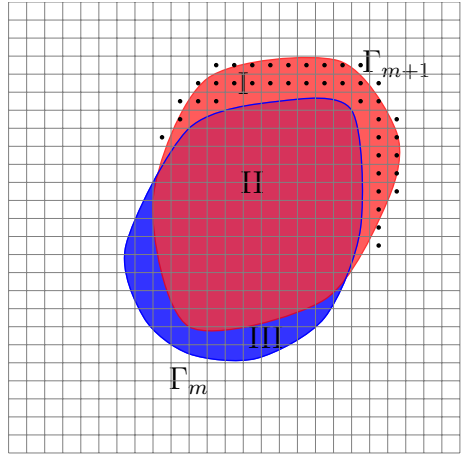

\noindent
{\bf Discretization of the reaction-diffusion equations on a moving cell region.}
Given $\Gamma_{m+1}$, $\Omega_{m+1}^+$, and $\Omega_{m+1}^-$, all 
specified by the level-set function $\phi^{m+1}$ on all the grid points, 
and also given the concentrations $u^m$ and $v^m$ on all the centers in $\Omega^+_m$ of grid cells, 
we need to find the approximate solution $u^{m+1}$ and $v^{m+1}$ 
on all the centers  of grid cells that overlap with $\Omega^+_{m+1}$ 
by the equations and boundary conditions \reff{eq:nonBC}. 
To do so, we first employ a second-order extrapolation method 
proposed in \cite{johansen1998cartesian} to extend $u^m$ and $v^m$ 
to the centers of grid cells that overlap with the  new cell region $\Omega^+_{m+1}$ but are 
not in $\Omega_m^+$; 
cf.\ black solid dotes in Figure~\ref{fig:cell_move}. 
We denote by $\tilde{u}_{i,j}^m$ and $\tilde{v}^m_{i,j}$ the original or an extended $u$-value and
$v$-value at the center of a grid cell in $\Omega_{m+1}^+$ labelled by $(i,j).$

Note by \reff{eq:nonEqu}--\reff{eq:nonBC} that 
the integral of $u+v$ over $\Omega$ is a constant with respect to time $t$  
(cf.\ \reff{M}), and its value is determined by the initial concentrations $u$ and $v$ at $t = 0.$
 We shall still denote this constant by $M$. 
To enforce this conservation of the total mass, 
we modify the value $\tilde{v}_{i,j}^m $ to get $v_{i,j}^m$ at 
centers of all the grid cells overlapping with $\Omega_{m+1}^+$ by 

\begin{align*}
v_{i,j}^{m} =\frac{1} {\text{Area}\,(\Omega_{m+1}^+)}
\left[ M - \int_{\Omega^+_{m+1}} \left( \tilde{u}^m + \tilde{v}^m \right) dA \right] 
+\tilde{v}_{i,j}^{m}.  
\end{align*}
The finally extended $u$ and $v$ values are now denoted by $u^{m}_{i,j}$ and $v^m_{i,j}$; they 
are defined on centers labelled by $(i,j)$ of grid cells overlapping with $\Omega_{m+1}^+$. 

We now focus on $u^{m+1}$ as $v^{m+1}$ is similar. We use the semi-implicit scheme 
\begin{equation*}
\frac{u^{m+1}-u^m}{\Delta t} =D_u \Delta u^{m+1}  + f(u^m, v^m) 
\qquad \mbox{for } (x, y) \in \Omega_{m+1}^+. 
\end{equation*}
Since the interface $\Gamma_{m+1}$ cuts through grid cells, we employ an 
embedded boundary method developed in \cite{papac2010efficient} to discretize the 
Laplacian operator.  Fix a grid cell $\mathcal{C}_{i,j}$ that overlaps with 
$\Omega_{m+1}^+$ and that is centered at $(x_i, y_j)$. 
Integrating both sides of the above equation over $\mathcal{C}_{i,j}\cap \Omega^+_{m+1}$, 
we have by the divergence theorem that 
\begin{equation}
\label{AreaLine}
\int_{\mathcal{C}_{i,j}\cap\Omega^+_{m+1}} \frac{u^{m+1}-u^m}{\Delta  t}\, dA
=D_u \int_{\partial(\mathcal{C}_{i,j}\cap \Omega^+_{m+1})} \nabla u^{m+1} \cdot \nu \, dl 
+\int_{\mathcal{C}_{i,j}\cap\Omega^+_{m+1}}f(u^m,v^m) \, dA,  
\end{equation}
where $\nu$ is the unit vector normal to the boundary $\partial (\mathcal{C}_{i,j} \cap \Omega^+_{m+1}).$ 
The two area integrals can be approximated by 
\begin{align}
\label{Area1}
&\int_{\mathcal{C}_{i,j}\cap\Omega^+_{m+1}}\frac{u^{m+1}-u^m}{\Delta t}\, dA
\approx \frac{u_{i,j}^{m+1}-u_{i,j}^m}{\Delta t}\, 
\mbox{Area}\,({\mathcal{C}_{i,j}\cap \Omega^+_{m+1}}), 
\\
\label{Area2}
& \int_{\mathcal{C}_{i,j}\cap\Omega^+_{m+1}}f(u^m,v^m)\, dA
\approx f(u_{i,j}^m, v_{i,j}^m) \, \mbox{Area}\,({\mathcal{C}_{i,j}\cap \Omega^+_{m+1}}). 
\end{align}
The area can be calculated using the level-set function $\phi^{m+1}$
 \cite{min2007geometric}. 
Whether or not the interface $\Gamma_{m+1}$ cuts through the grid cell $\mathcal{C}_{i,j}$, 
by the boundary condition $\partial_n u^{m+1} = 0$ on $\Gamma_{m+1}$, 
we can approximate the line integral in \reff{AreaLine} by \cite{papac2010efficient}
\begin{align}
\int_{\partial(\mathcal{C}_{i,j}\cap \Omega^+_{m+1})} \nabla u^{m+1} \cdot \nu \, dl 
& \approx
\frac{u_{i+1,j}^{m+1}-u_{i,j}^{m+1}}{h} L_{i+1/2,j}
-\frac{u_{i,j}^{m+1}-u_{i-1,j}^{m+1}}{h} L_{i-1/2,j}
\nonumber 
\\
\label{LineIntegral}
& 
+\frac{u_{i,j+1}^{m+1}-u_{i,j}^{m+1}}{h} L_{i,j+1/2 }
-\frac{u_{i,j}^{m+1}-u_{i,j-1}^{m+1}}{h} L_{i,j-1/2 }, 
\end{align}
where $L_{i\pm 1/2, j} \in [0,h]$ and $L_{i, j\pm 1/2} \in [0,h]$ refer to the length 
of the corresponding edge of the grid cell $\mathcal{C}_{i,j}$ inside $\Omega^+_{m+1}$. 
These lengths can be calculated using the level-set function $\phi^{m+1}$ that defines
the interface $\Gamma_{m+1}$ \cite{min2007geometric}. 

The coefficient matrix of the resulting system of linear equations is 
symmetric positive definite \cite{papac2010efficient}, and the system can be solved by 
the conjugate gradient method with an incomplete Cholesky preconditioner or by an algebraic multigrid method.

\medskip

\noindent
{\bf Algorithm.}
\begin{compactenum}
\item[Step 0.]
Input all the parameters. 
Set the computational box $\Omega = (-L, L)^2 \subset \R^2$ and cover it with a uniform
finite-difference grid with grid sizes $h$. 
Discretize the time interval $[0, T]$ of interest with time step $\Delta t.$ 
Initialize the level-set function $\phi^0$ and the concentrations $u^0$ and $v^0.$ 
Set $m = 0.$ 
	
\item[Step 1.]
Extend the normal velocity from the interface to the entire computational box.
Solve the semi-implicit discretization equation \reff{phiphi} to 
get the updated level-set function $\phi^{m+1}$. 
Reinitialize the level-set function and still denote it by $\phi^{m+1}$. 
	
\item[Step 2.]
Extend $u^m$ and $v^m$ to centers of grid cells overlapping with $\Omega_{m+1}^+$ defined by $\phi^{m+1}.$ 
Solving the semi-implicit discretization equations (cf.\ \reff{AreaLine}--\reff{LineIntegral})
to obtain $u^{m+1}$ and $v^{m+1}$. 

\item[Step 3.]
Check if the cell region $\Omega_{m+1}^+$ touches the boundary $\partial \Omega$. 
If so, shift the computational box so that the cell
is centered in the new computational box, still denoted $\Omega$. 
	
\item[Step 4.]
Set $m:=m+1$. Repeat Steps 1--3 until the final simulation time is reached.  
	
\end{compactenum}


\medskip

\begin{figure}[h]
\centering
\includegraphics[width=0.5\linewidth]{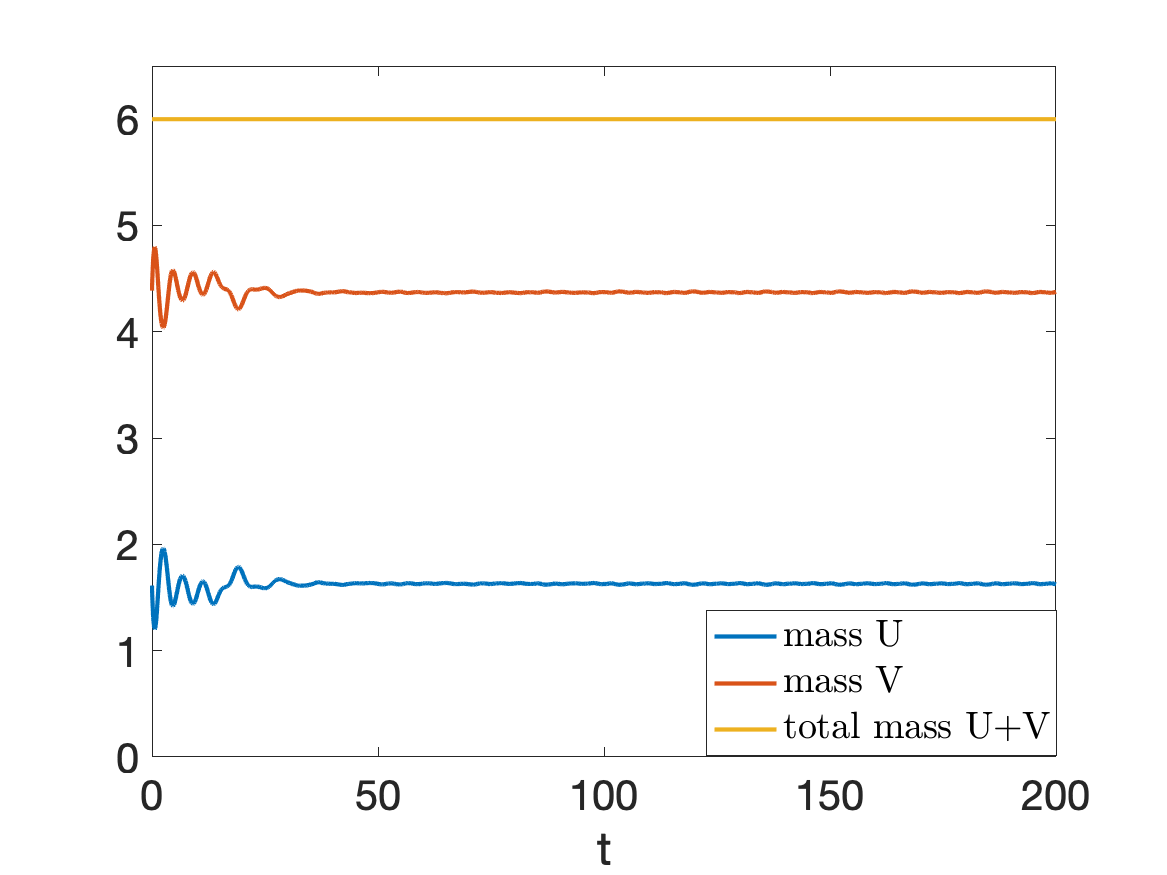}
\caption{Numerical conservation of total mass. Here $U = U(t)$ and $V= V(t)$ at time $t$ are 
defined as the integral of $u(x, y, t)$ and $v(x, y, t),$ respectively, over the cell region 
$\Omega^+(t)$, and the total mass is defined to be the sum $U + V$.}
\label{f:mass}
\end{figure}

\noindent 
{\bf Convergence test.} 
We have tested our numerical methods and code. In Figure~\ref{f:mass}, we show that the total
mass conservation is captured numerically in a long time simulation. 

\begin{figure}[h]
	\centering
\includegraphics[width=0.47\linewidth]{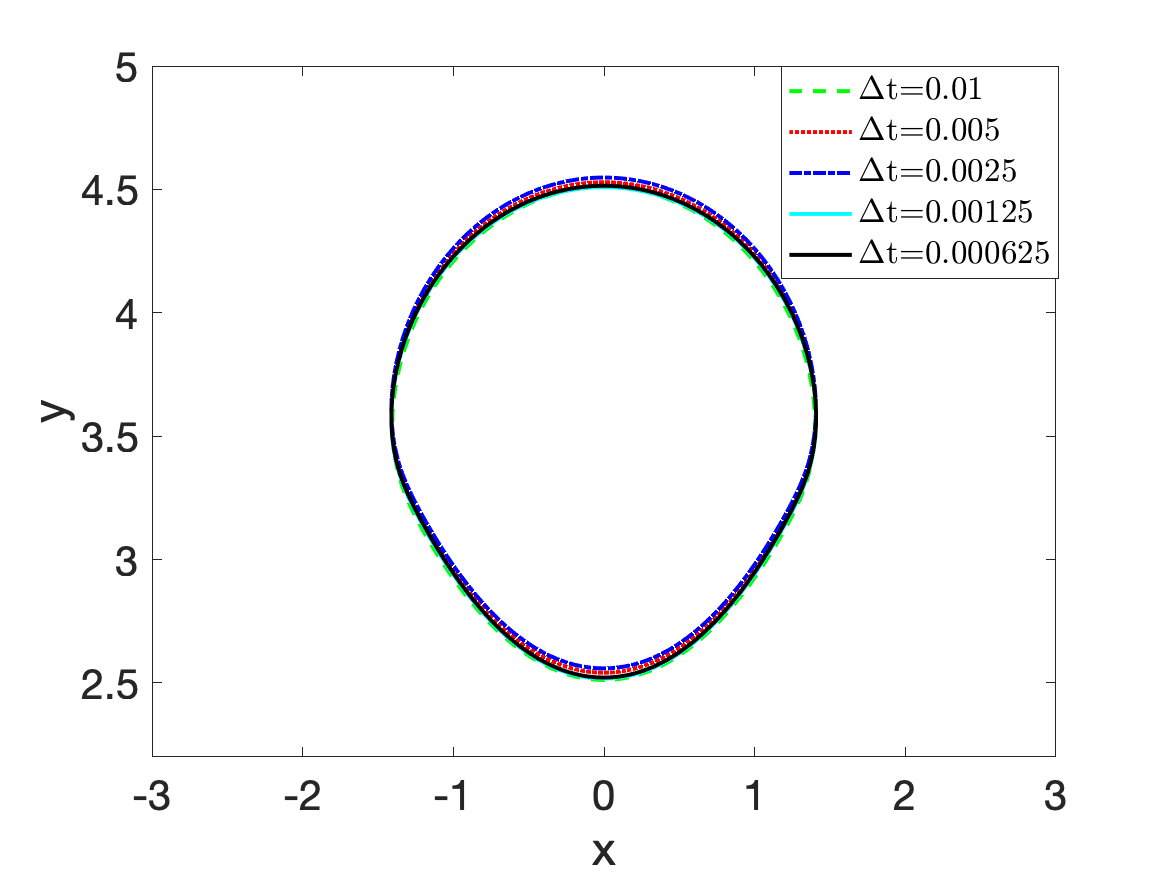}
\includegraphics[width=0.47\linewidth]{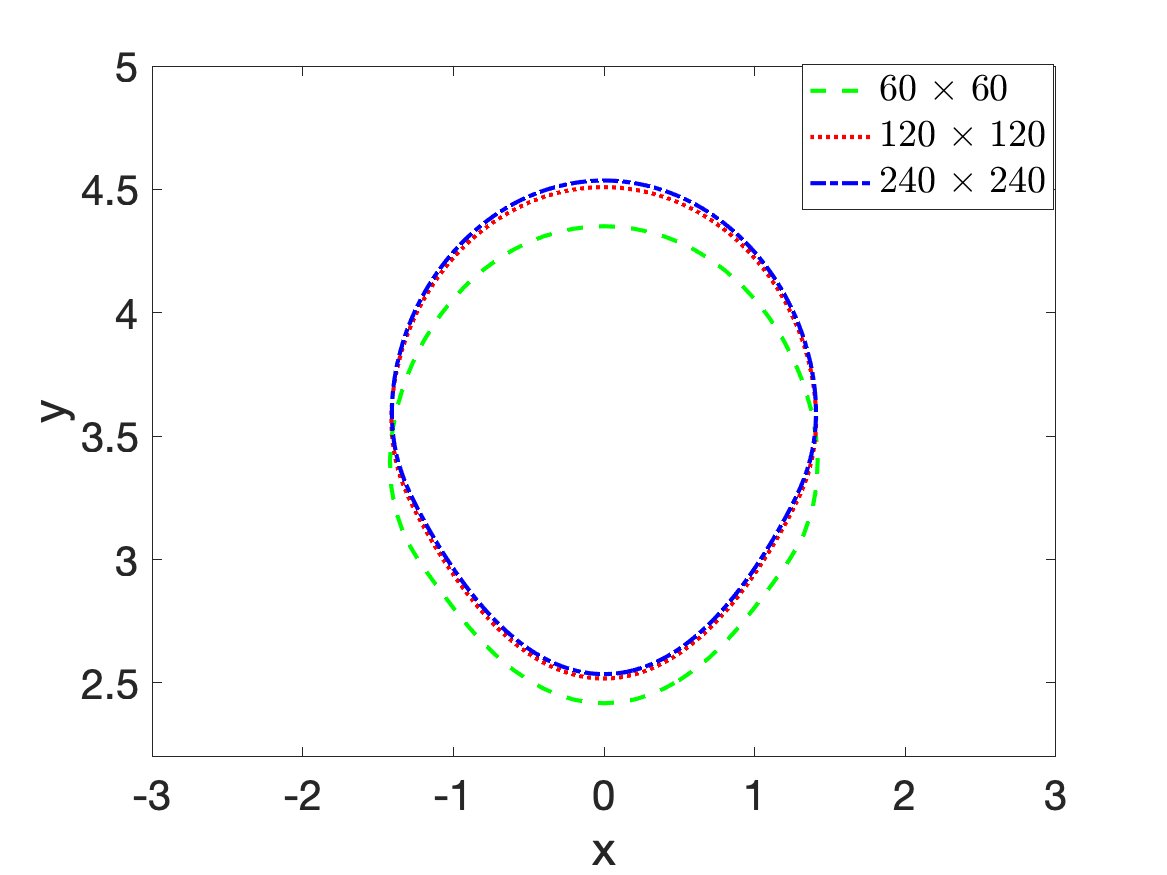}
	\caption{Simulated cell shapes with 
the final simulation otime $T = 10.$  
The initial simulation box $\Omega = (-3, 3)^2$ is shifted during the simulation. 
Left:  Simulation results with the grid size $h= 0.05$ and different time steps $\Delta t$.  
Right: Simulation results with the time step  $\Delta t = 1.25\times10^{-3}$ and with different 
number of grid points.}
	\label{fig:converge_study}
\end{figure}

We have also used our numerical methods to 
simulate a moving cell with the final (rescaled) time being $T = 10$, 
using  different time steps and different spatial grid sizes. 
Figure~\ref{fig:converge_study} shows our simulation results. They indicate that 
our numerical method and algorithm converge both in time and space.

%% file: SimulationResults.tex
\section{Simulation Results and Analysis}
\label{s:Results}

We perform numerical simulations to study the cell polarization in response to various external 
stimulus and the  trajectory of a moving cell, and analyze these simulation results 
in terms of the modelling and parameters.


\subsection{Parameters and Non-dimensionalization}

In Table~\ref{parameters}, we collect all the parameters
in the original model \reff{NormalVelocity}--\reff{BC4uv}, 
and describe their meanings and units. We also provide 
their estimated values following \cite{camley2017crawling, ShaoRappelLevine_PRL10}.

\begin{table}[h]
\caption{Parameters}\label{parameters}
\vspace{-5mm}
\begin{center}
\scalebox{1}{
\begin{tabular}{ l l l l}
\hline
Parameters & Description &Estimated Values& Units\\
\hline

$D_u$ &  diffusion coefficient of  $u$  & $0.1 \sim 0.5 $ & $ \mu {\rm m}^2/{\rm s}$ \\
$D_v$ &  diffusion coefficient of  $v$ & $ 10 \sim 50 $ & $\mu {\rm m}^2/{\rm s}$ \\

$\alpha$ &coefficient of F-actin extension&  0.1&$ {\rm pN}/\mu {\rm m}$\\

$\beta$&coefficient for myosin retraction& 0.2 &${\rm pN}/\mu {\rm m}$ \\

$\tau$ & friction coefficient &2.62&$ {\rm pN}{\rm s}/\mu {\rm m}^2$\\

$\gamma$ & surface tension &1&${\rm pN}$\\
$k$ & relative reaction rate &  $\sim \, $0.01&${\rm s}^{-1}$\\
$c$ & concentration of  $u$ at the cell front& 1 $\sim$ 10 &concentration unit\\
$C$  & interconversion  parameter &0.5 $\sim$ 0.8&unitless\\
\hline
\end{tabular}
}
\end{center}
\end{table}

To non-dimensionalize our equations \reff{NormalVelocity}--\reff{BC4uv}, we 
  follow \cite{camley2017crawling} to introduce two parameters. 
One is the typical cell speed $V_0$ which is in the range $\sim\,$0.1$ \mu {\rm m}/s$. 
The other is the typical radius of a cell $R$ which is in the range 
$\sim\,$10$\mu {\rm m}.$ We then 
introduce non-dimensionalized parameters according to Table \ref{unitless}. 
\begin{table}[h]
\caption{Nondimensionalized  Parameters}\label{unitless}

\vspace{-5mm}

\begin{center}
\begin{tabular}{ l l c}
\hline
Parameters & Description &Estimated Values\\
\hline

\vspace{1mm}

$\widehat{D}_u=\frac{D_u}{V_0{R}}$ &  Rescaled diffusive coefficient of  $u$ & $ 0.1 \sim 0.5 $ \\

\vspace{1mm}

$\widehat{D}_v=\frac{D_v}{V_0{R}}$ &  Rescaled  diffusive coefficient of  $v$ & $ 10 \sim 50 $\\

\vspace{1mm}

$ K=\frac{k{R}\rm c^2}{V_0}$ & Rescaled reaction rate  compared to motility & $ 100 \sim 500 $ \\

\vspace{1mm}

$\chi =\frac{\gamma}{V_0\tau {R}}$ & Relative strength of surface tension & $ 0.1 \sim 0.3$ \\

\vspace{1mm}

$\widehat{u^*}=\frac{\beta}{{\rm c}\alpha}$ & Rescaled contractility & $ 0.2 \sim 0.45$ \\

\vspace{1mm}

$\widehat{M} =\frac{N_{tot}}{\rm c{R}^2}$ & Rescaled total amount of protein $u$ and $v$ & $6 \sim 8$ \\
\hline
\hline
\end{tabular}

\end{center}
\end{table}

We now define $\hat{x}=x/R$, $\hat{t} = (V_0/R) t$,   $\hat{u}=u/c $, $\hat{v}=v/c$, 
and $\widehat{V}=V/V_0$, and convert the original system of equations
\reff{NormalVelocity}--\reff{BC4uv} into the non-dimensionalized system of 
equations for $\widehat{V}$, $\hat{u}$, and $\hat{v}$, which is, after dropping all the hat, 
the system \reff{eq:nonV}--\reff{eq:nonBC}.

\subsection{Cell Polarization}
\label{ss:CellPolarization}

Inspired by the one-dimensional simulations of the wave-pinning mechanism \cite{mori2008wave}, 
we consider a non-moving or stationary cell that occupies the fixed region 
$\Omega^+ $ whose boundary is the curve 
\[
x=(1-0.3\cos2\theta)\cos\theta, \quad  
	y=(1-0.3\cos2\theta)\sin\theta \qquad \forall \, \theta\in[0, 2\pi).
\]
We also take the total amount of proteins (i.e., the total mass of the two species) 
to be $M = 6$ after nondimensionalization. The computational box is $\Omega = (-2.5, 2.5)\times (-2.5, 2.5)$.  Note that $\Omega^+ \subset \Omega.$ 
We cover $\Omega$ by uniform finite-difference grid of grid size $h=0.05$.  Other parameters are taken as $\triangle t=0.001$, $D_u=0.3,\, D_v=30$, $K=500$, and $C=0.8$.

\medskip

\noindent
{\bf Random initial value.}
We choose the initial value $u_0 = u(\cdot,0)$ to be a random variable defined on all the 
grid centers. The values are generated uniformly at random from $[0, 0.8].$
We also set the initial value $v_0 = v(\cdot, 0)$ to be a constant, which is 1.342 in our case with total mass $M=6$. We then solve numerically the reaction-diffusion equations with the zero Neumann boundary conditions. 

\begin{figure}[h]
	\centering
	\includegraphics[width=0.323\textwidth]{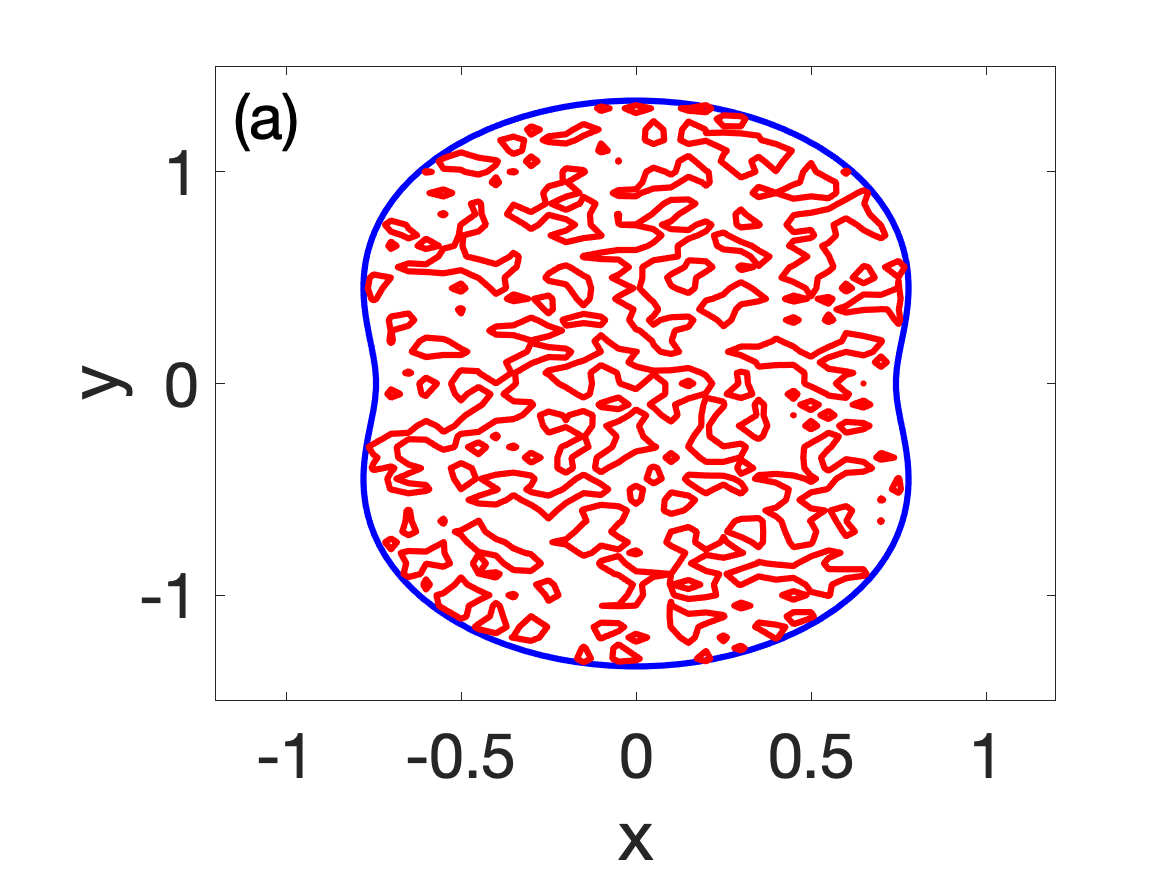}
	\includegraphics[width=0.323\textwidth]{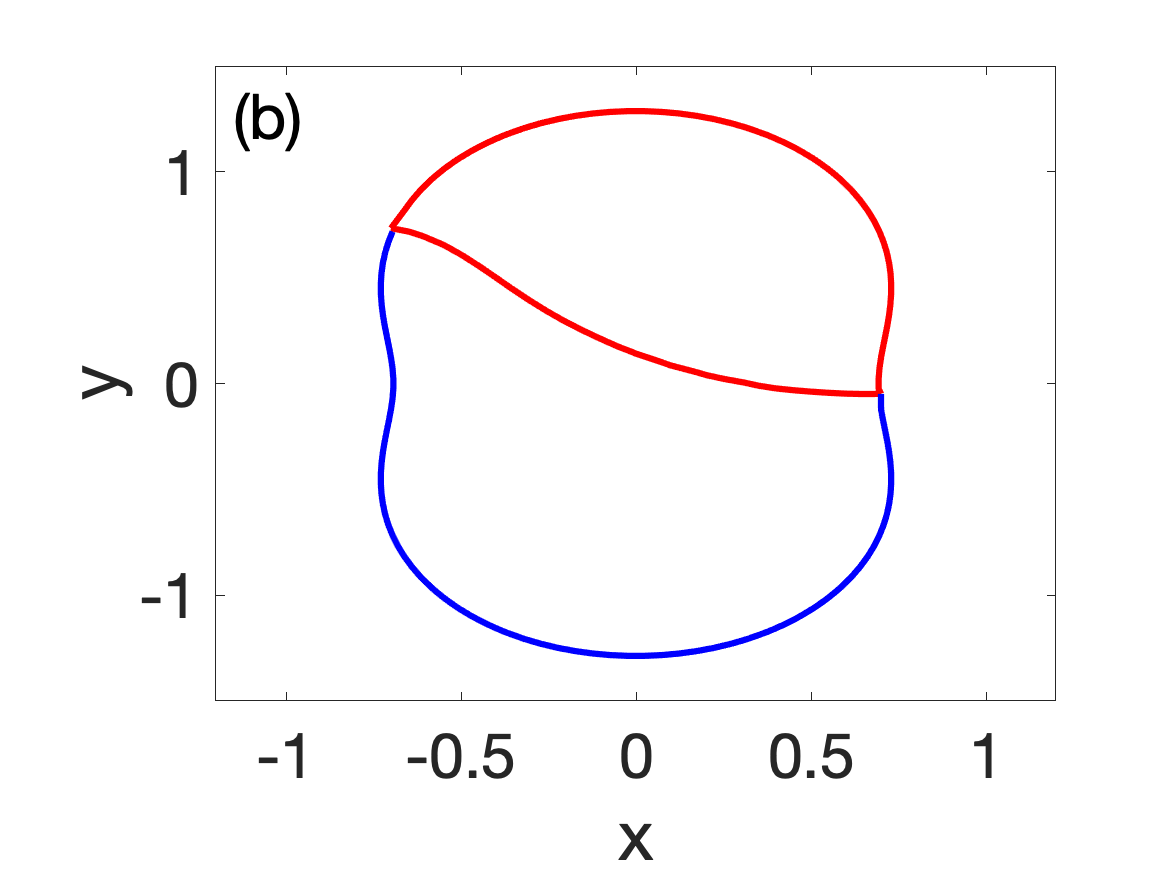}
	\includegraphics[width=0.323\textwidth]{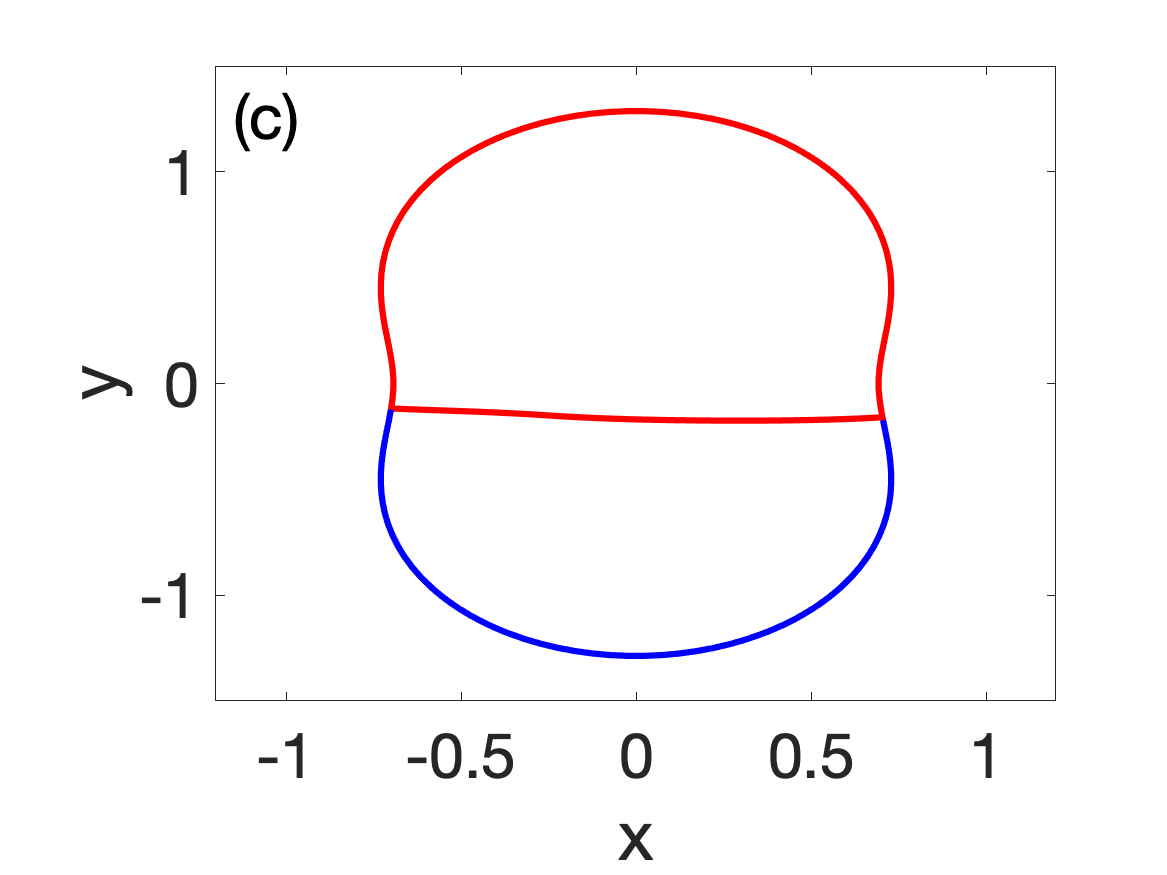}
	\caption{Cell polarization with a random initial value $u_0.$ 
		(a) Cell region bounded by the blue curve. A random initial concentration field $u_0$ 
 is distributed at grid points inside the cell region. 
Red spot is the set of points at which $u_0 \ge 0.5.$ 
		(b) At $t = 0.5$, the region of $u \ge 0.5$, 
bounded by the red curves, is located in the front of the cell, showing a polarized cell. 
		(c) At $t = 2$, the cell polarization has reached an equilibrium.
	}
	\label{fig:random_noise}
\end{figure}

Figure~\ref{fig:random_noise} (a) shows the cell region bounded by the blue curve. 
The red region inside the cell is the set of points at which the initial random concentration
value $u_0 \ge 0.5.$

Figure~\ref{fig:random_noise} (b) shows the part of the cell region, marked by the closed red curve,
at which the concentration $u \ge 0.5.$ This shows that the cell is polarized at this (rescaled) time
$t = 0.5$. 

As discussed in \cite{jilkine2009wave,camley2017crawling},  
the reaction-diffusion system tends to minimize the length of the interface that separates the high and low $u$-concentration regions.

Figure~\ref{fig:random_noise} (c) shows that at 
$t=2$ the interface between the high and low $u$-concentration regions 
does not change, indicating the cell polarization reaches an equilibrium \cite{mori2008wave}. 

\medskip

\noindent
{\bf External stimulus.}
We introduce an external stimulus and solve the system of equations
\begin{align*}
	&\partial_t u =D_u  \Delta u+f(u,v)+ S v \quad \mbox{in } \Omega^+ \times (0, T], 
	\\
	&\partial_t v =D_v  \Delta v-f(u,v)- S v \quad \mbox{in }  \Omega^+ \times (0, T], 
\end{align*}
with the same boundary conditions $\partial_n u = \partial_n v = 0$ on $\partial \Omega^+$ 
and a final simulation time $T$. 
The stimulus function is $S v$ with $S$ defined on $\Omega^+ \times [0, T]$ by 
\begin{equation*}
	S = S (x, y, t) 
	= \begin{cases}
		s_1(t)(1.3-y)(0.7-x) \quad &\text{if} \, (x,y,t )\in \Omega^+ \times [0, 1], \\
		s_2(t)(y+1.3)(x+0.7) \quad &\text{if} \, (x,y,t)\in \Omega^+ \times [10, 11], \\
		0 \quad &\text{elsewhere}, 
	\end{cases}
\end{equation*}
where 
\begin{equation*}
	s_1(t)=\begin{cases}
		0.07\quad &\text{if} \ 0 \le t \le 0.5, \\
		0.07(1-\frac{t-0.5}{0.5}) \quad  &\text{if} \ 0.5 < t \leq  1, \\
		0 \quad &\text{elsewhere,}
	\end{cases}
\quad \mbox{and}  \quad 
s_2(t) = \begin{cases}
  s_1(t-10) \quad & \text{if} \  10 \le t \le 11, \\
  0 \quad & \text{elsewhere.} 
\end{cases}
\end{equation*}
Note that the stimulus is strong in certain region of the cell
during the time interval $0 \le t \le 1$ and 
is strong in a different region of the cell during $10 \le t \le 11.$ 

\begin{figure}[h]
	\centering
	\includegraphics[width=0.3\textwidth]{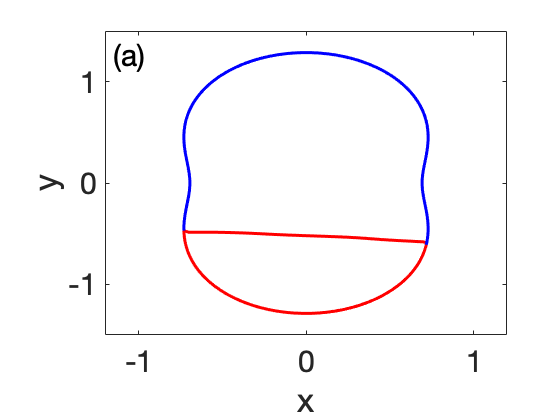}
	\includegraphics[width=0.3\textwidth]{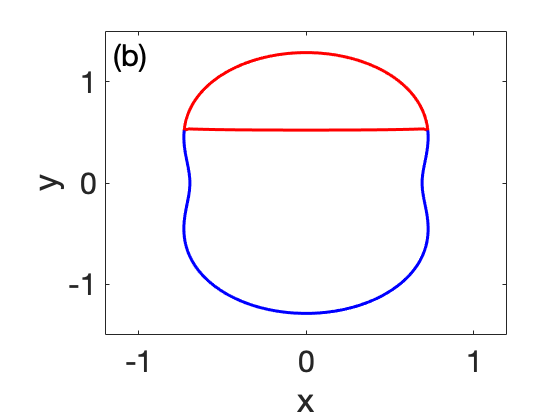}
	\caption{
A stimulus, strong in a local region of the cell, 
 is applied to the initially homogeneous state
during $0 \le t \le 1$ and then turned off during $ 1 < t < 10$. 
The stimulus, strong in a different region of the cell,  is later turned on 
during $10 \le t \le 11$ and turned off again after $t = 11$ till the final simulation time $T  = 20.$ 
(a) At $t = 10$, the cell is polarized. 
(b)  At $t = 20$, the cell is polarized again but with a reversed orientation.}
\label{fig:graded_stimuli}
\end{figure}

We set the initial values of $u$ and $v$ to be constant and solve the reaction-diffusion
equations with the stimulus up to the final simulation time $T = 20.$  
When the locally strong stimulus is turned on from $t=0$ to $t=1$,
the active form $u$ increases locally near the south-west corner. 
Such increase then leads to the formation of a spatial interface in the cell, 
separating the high and low concentrations of the active form $u$, 
that propagates inside the cell region. 
Meanwhile, the concentration $v$ of the inactive form decreases, 
leading accordingly to the decreasing of $u^+ = C v$ in the kinetic form of $f$.
As a result, the motion of the internal interface slows down, and 
is finally pinned down, and the cell reaches a polarized steady state \cite{mori2008wave}. 
Figure~\ref{fig:graded_stimuli} (a) shows such a polarized state at $t = 10.$ 
The cell is re-polarized, with a reversed orientation,
 after the stimulus is turned on again during $10 \le t \le 11$
but strong in a different spatial region of the cell; cf.\ Figure~\ref{fig:graded_stimuli} (b). 


\subsection{Cell Trajectory }
\label{ss:CellTrajectory}

We define the cell trajectory of a moving cell 
to be the time trajectory of the geometrical center of the cell $(X_{\rm c}(t), Y_{\rm c}(t))$, which is defined by 
\begin{equation*}
X_{\rm c}(t)= \frac{1}{\mbox{Area}\,(\Omega^+(t)) } \int_{\Omega^+(t)} x \, dx dy 
\quad \mbox{and} \quad 
Y_{\rm c}(t)= \frac{1}{\mbox{Area}\,(\Omega^+(t)) } \int_{\Omega^+(t)} y \, dx dy.  
\end{equation*} 
We study two typical types of trajectories, straight 
 and circular trajectories, aiming at a qualitative understanding of controlling parameters for such trajectories. 

We shall also compare the two-species and one-species models in terms of
the prediction of different trajectories. 

 In all the simulations reported below, the initial cell region  is a circle of radius
 $1.3$, which is centered at $(0,-1)$, with a radius $1.3$. The cell is polarized with the concentration of $u$ to be 0.8 where $y \geq -0.8$, and $u=0$ in the remaining part, while $v$ is uniformly distributed on the cell domain, with the total mass $M=6$.   We also set  the kinetic parameters $K=100$ and  $C=0.8$, the grid size 
$h=0.06$ and time step $\Delta t=0.005.$

\medskip

\begin{figure}[h]
	\centering
\includegraphics[width=0.4\linewidth]{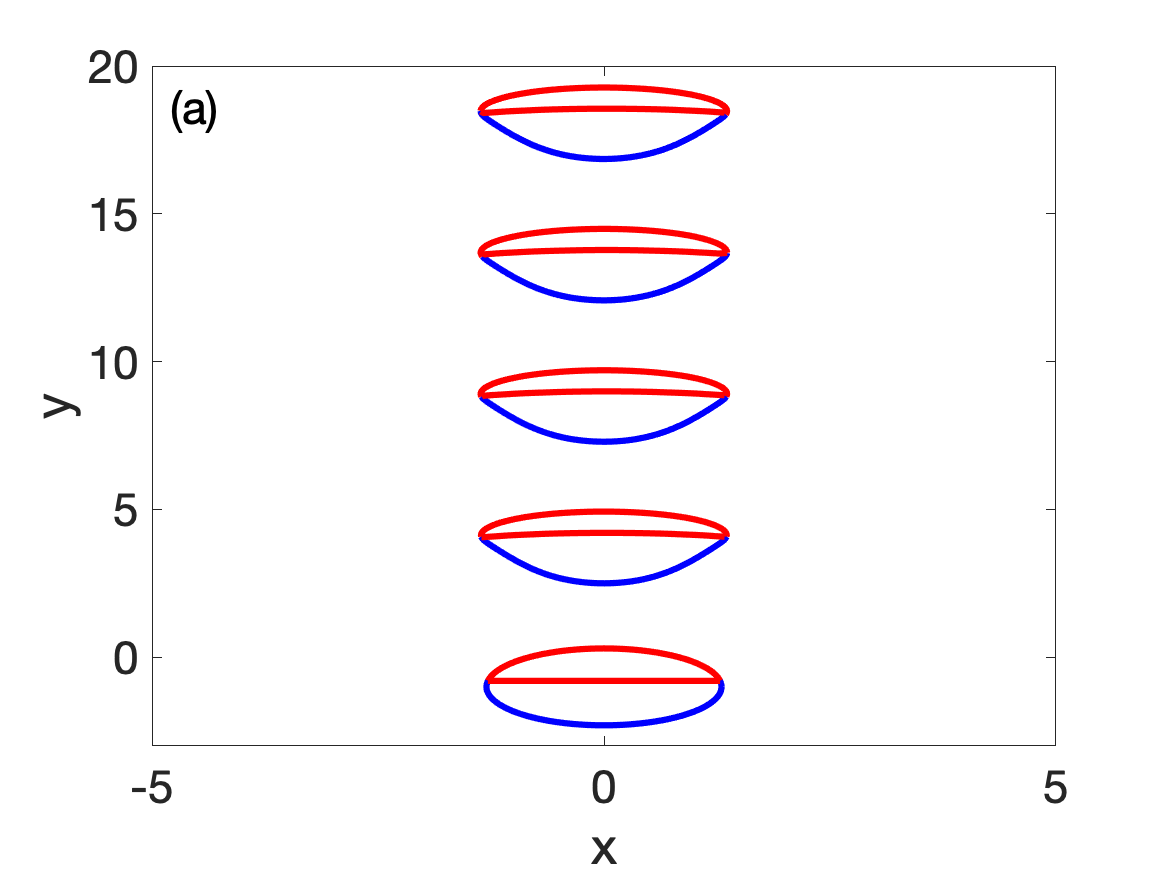}
\includegraphics[width=0.4\linewidth]{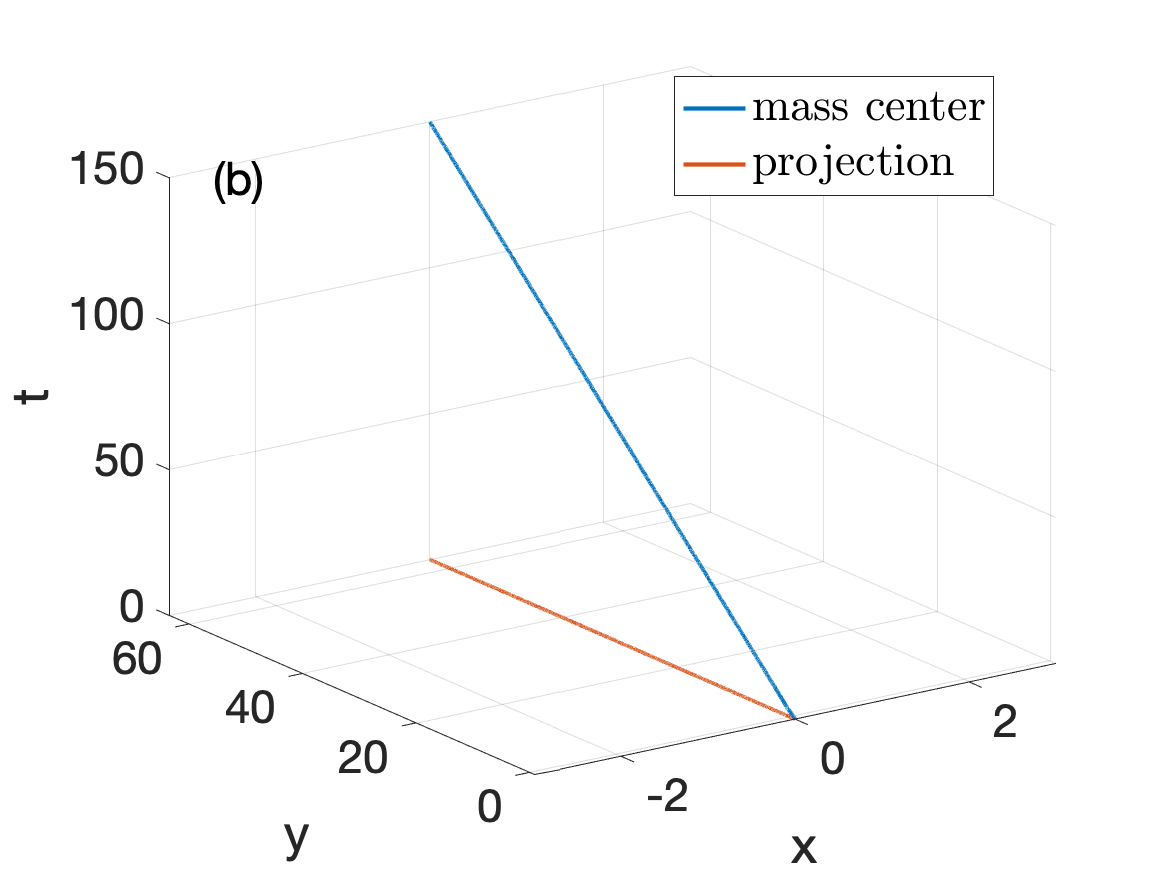}
\includegraphics[width=0.4\linewidth]{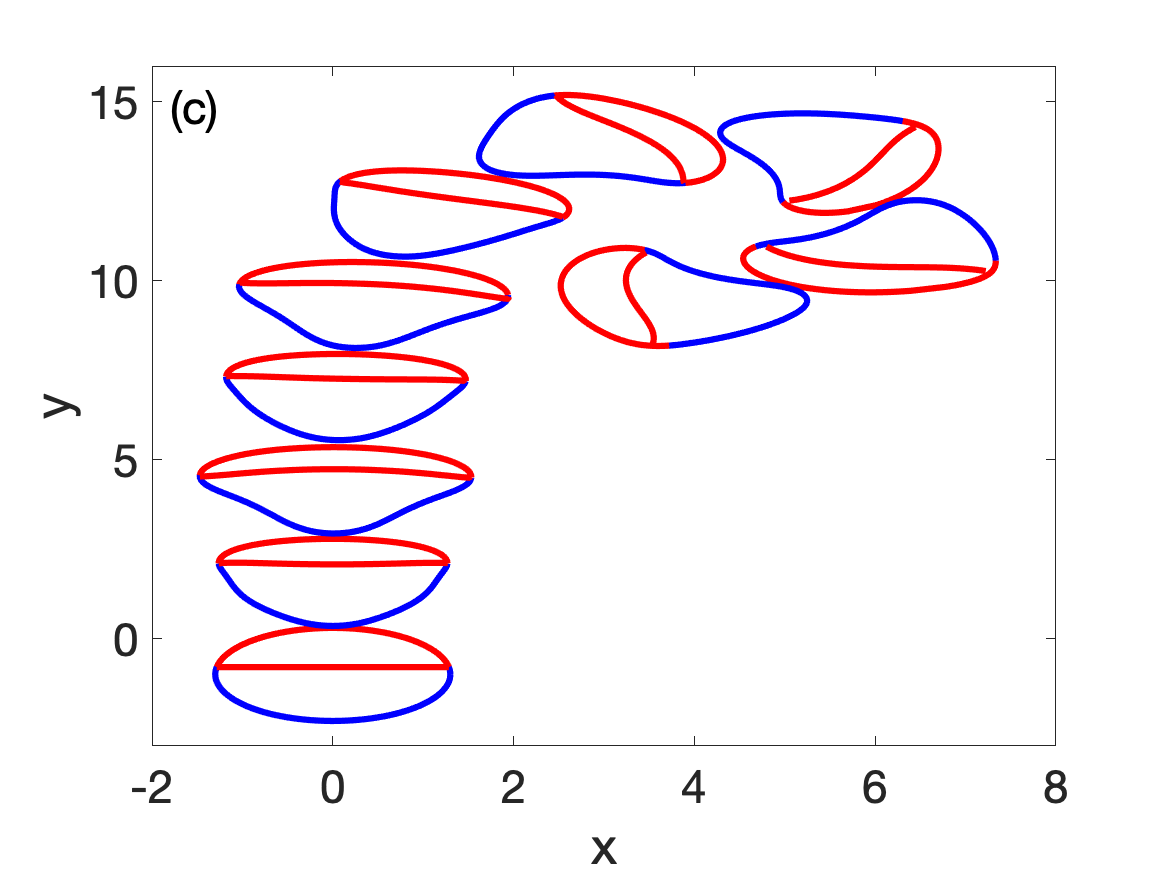}
\includegraphics[width=0.4\linewidth]{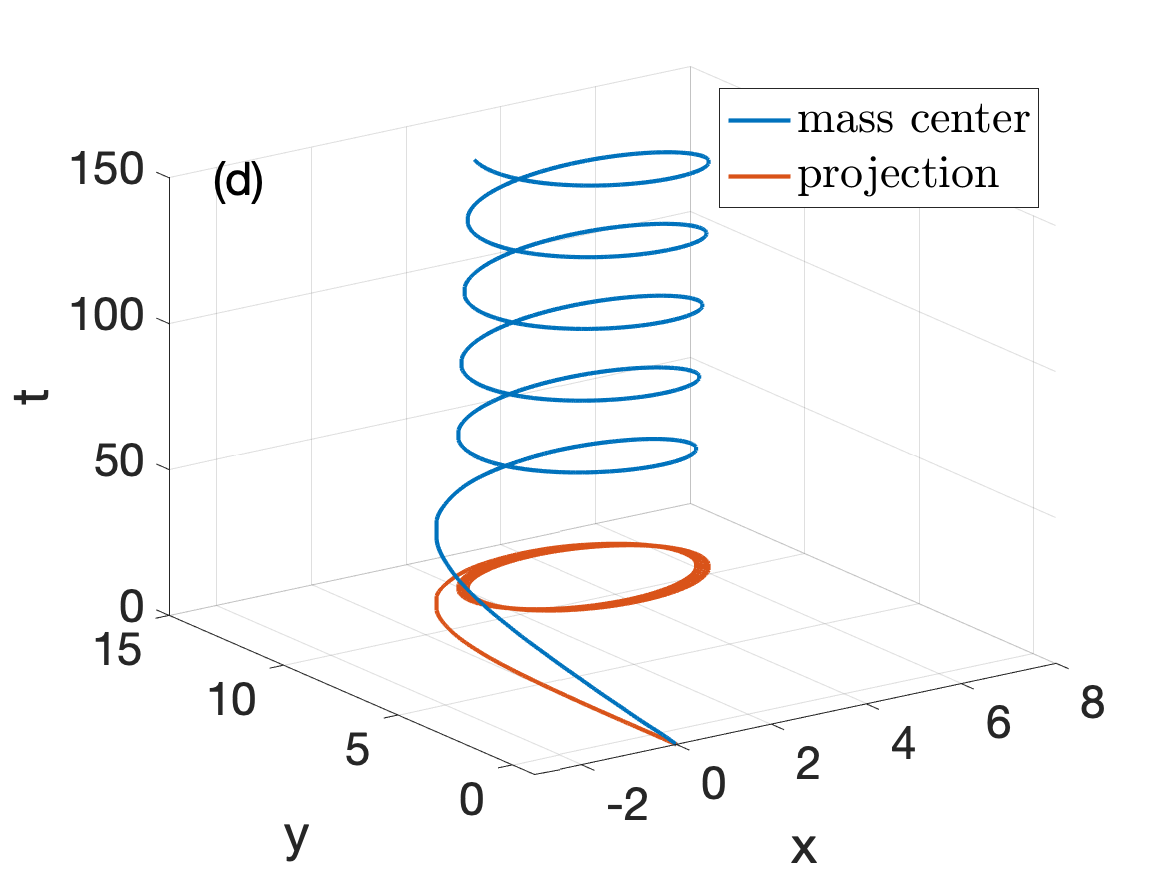}
	\caption{Long-time cell movement simulations. 
In (a) and (b), $D_u=0.1,$ $D_v=10,$ $\chi=0.2, $ and  $u^*=0.2$.  
In (c) and (d), $D_u=0.5,$ $D_v=50,$ $ \chi=0.1,$ and $u^*=0.25$. 
The sequence of snapshots of cells in (a)  are taken at 
$t=0, \,10, \,20,\,30, \,40$ and those in (c) are taken at 
$t=0, \, 5, \, 10, \, 15,  \,20, \, 25, \,30, \, 35, \,40, \, 45$. 
The blue line or curve in (b) or (d) is the space-time cell trajectory, 
while the red line or curve in (b) or (d), marked ``projection" is the (two-dimensional)
space trajectory of the moving cell. 
}
	\label{fig:celltrajectories}
\end{figure}

\noindent
\textbf{Long-time trajectories.}
We simulate a moving cell with two different sets of parameters $D_u, $ $D_v$, and $\chi$, and 
plot the cell trajectory in Figure~\ref{fig:celltrajectories}. 
We observe clearly a straight trajectory 
(cf.\  Figure~\ref{fig:celltrajectories} (a) and (b))
and the circular trajectory (cf.\  Figure~\ref{fig:celltrajectories} (c) and (d)). 
Note that the parameters we use in these simulations 
are similar to those used in \cite{camley2017crawling} to capture 
both the straight and circular trajectories as in Figure \ref{fig:celltrajectories} with the final time $t = 40$ while here we have simulated the cell movement up to $t = 150,$ 
indicating that the two patterns are persistent, and the model and our methods are robust.


\medskip

\begin{figure}[h]
	\centering
\includegraphics[width=0.51\linewidth]{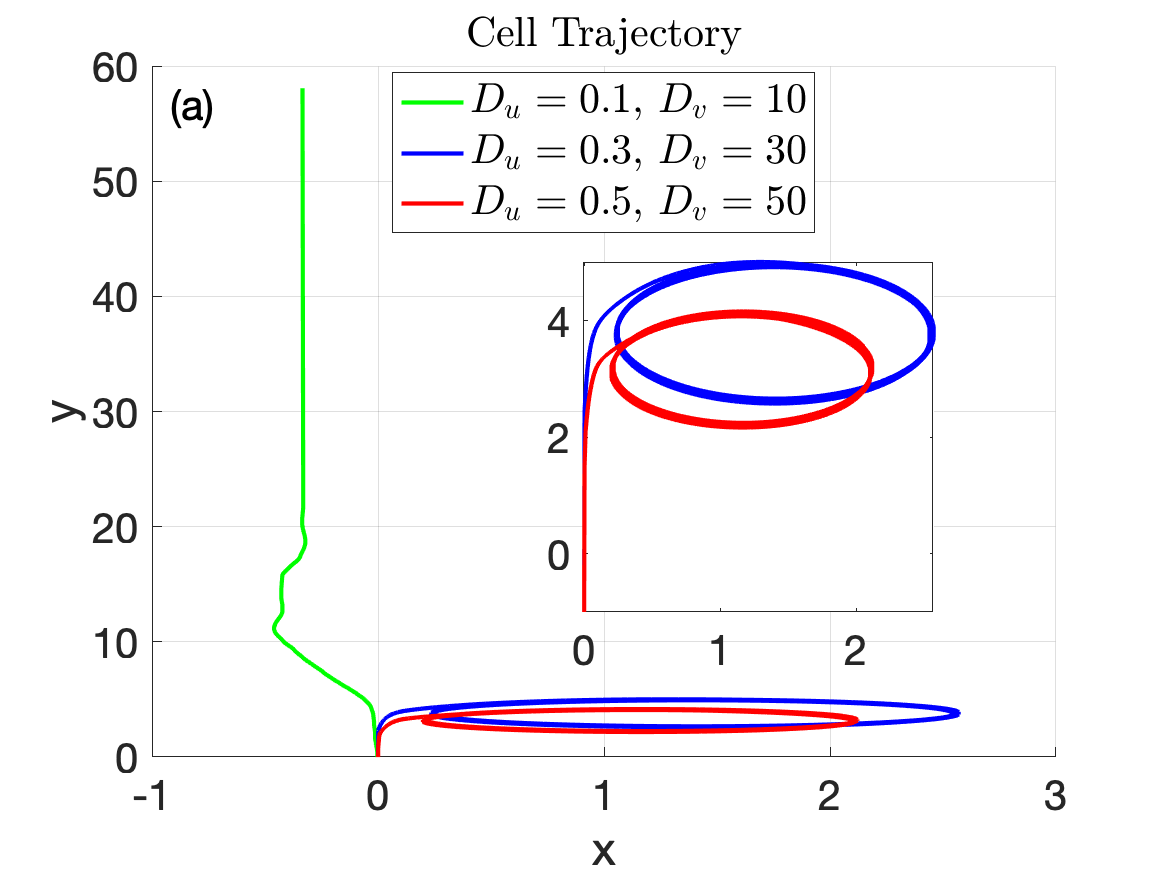}
\hspace{-20 pt}
\includegraphics[width=0.51\linewidth]{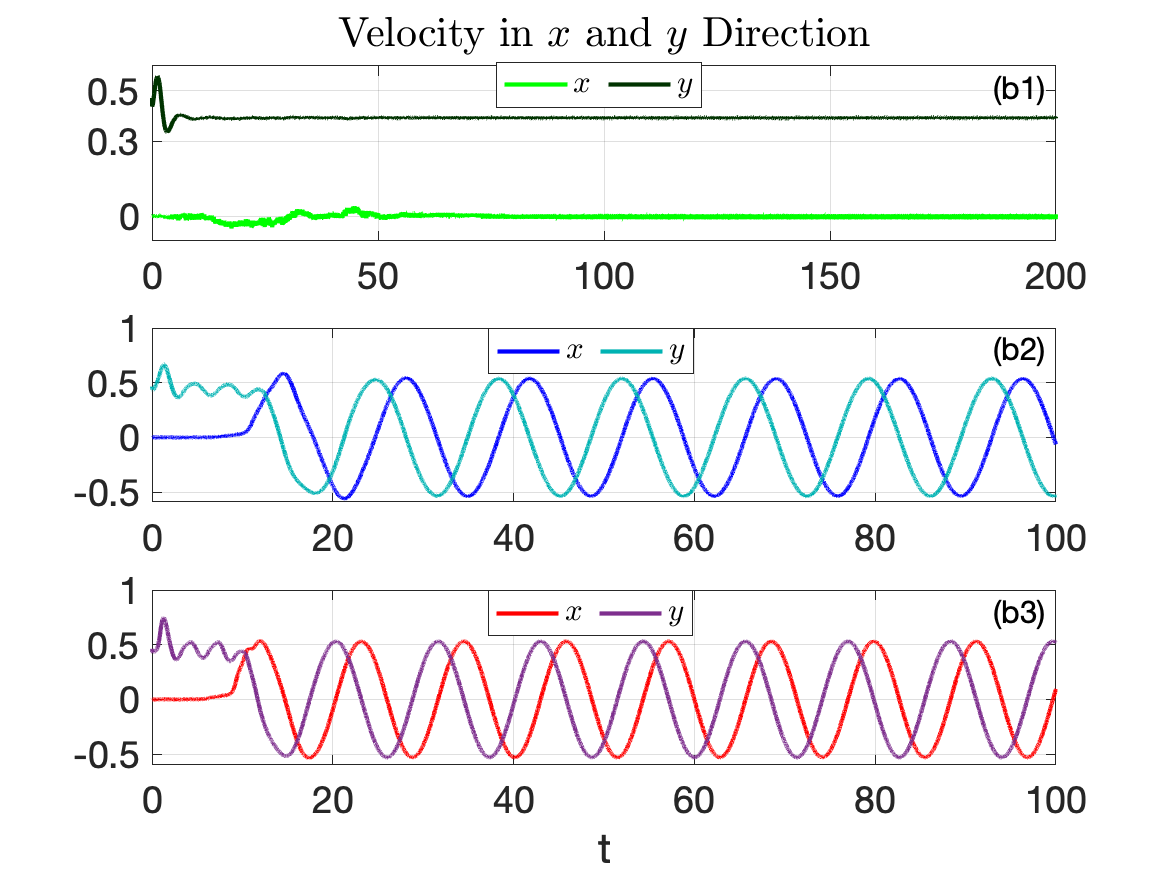}
	\caption{
(a) Cell trajectories predicted with three different sets of diffusion constants $D_u$ and $D_v.$ 
The small window is the zoom-in of the two circular trajectories. 
(b1)--(b3) The $x$ and $y$ components of the velocity at the center of a moving cell predicted by 
our numerical simulations corresponding to the three sets of $D_u$ and $D_v$ values marked in (a). 
}
	\label{fig:diffusionrate}
\end{figure}

\noindent
{\bf Effects of diffusion.}
We now set $\chi=0.1$ and $u^*=0.4$, and vary the diffusion constants
$D_u$ and $D_v$ to study how the diffusion can affect the cell movement. 
In Figure~\ref{fig:diffusionrate}, we plot our simulation results for three sets of 
diffusion constants: 
case 1:  $D_u=0.1$ and $ D_v=10$; 
case 2: $D_u=0.3$ and $ D_v=30$; 
case 3: $D_u=0.5$ and $ D_v=50$.  
In Figure~\ref{fig:diffusionrate} (a), we observe that 
the cell trajectory is linear (or straight) for case 1, while it is 
circular for case 2 and case 3. 
Note that the plot in the window
is the zoom-in of the long-time trajectories for case 2 and case 3.  
In all the three cases, there is a preparation time before the cell starts to move 
in a straight line for case 1 or in a circular pattern for case 2 and case 3. 

In Figure~\ref{fig:diffusionrate} (b1)--(b3), we plot the $x$ and $y$ components of 
the velocity at the geometrical center of the moving cell corresponding to the three cases, 
respectively. We observe that larger diffusion constants correspond to a shorter preparation time before
the onset of the linear or circular trajectory. Moreover, fast diffusion is correlated to 
a smaller circular trajectory. 


\medskip

\noindent
{\bf Contractility.}
This refers to the cell contraction due to the decreasing of concentration $u$ 
the rear part of the cell. In the model, the cell contractility is determined by the 
threshold concentration $u^\ast$. 
To study how the variation of $u^\ast$ can affect the cell trajectory, we fix the 
diffusion constants $D_u=0.4$ and $D_v=40$ and the rescaled surface tension constant $\chi=0.1$,  
and simulate the cell movement with different values of $u^\ast$: 
$0.25,$ $0.3$, and $0.4.$ 

\begin{figure}[H]
	\centering
		\includegraphics[width=0.5\linewidth]{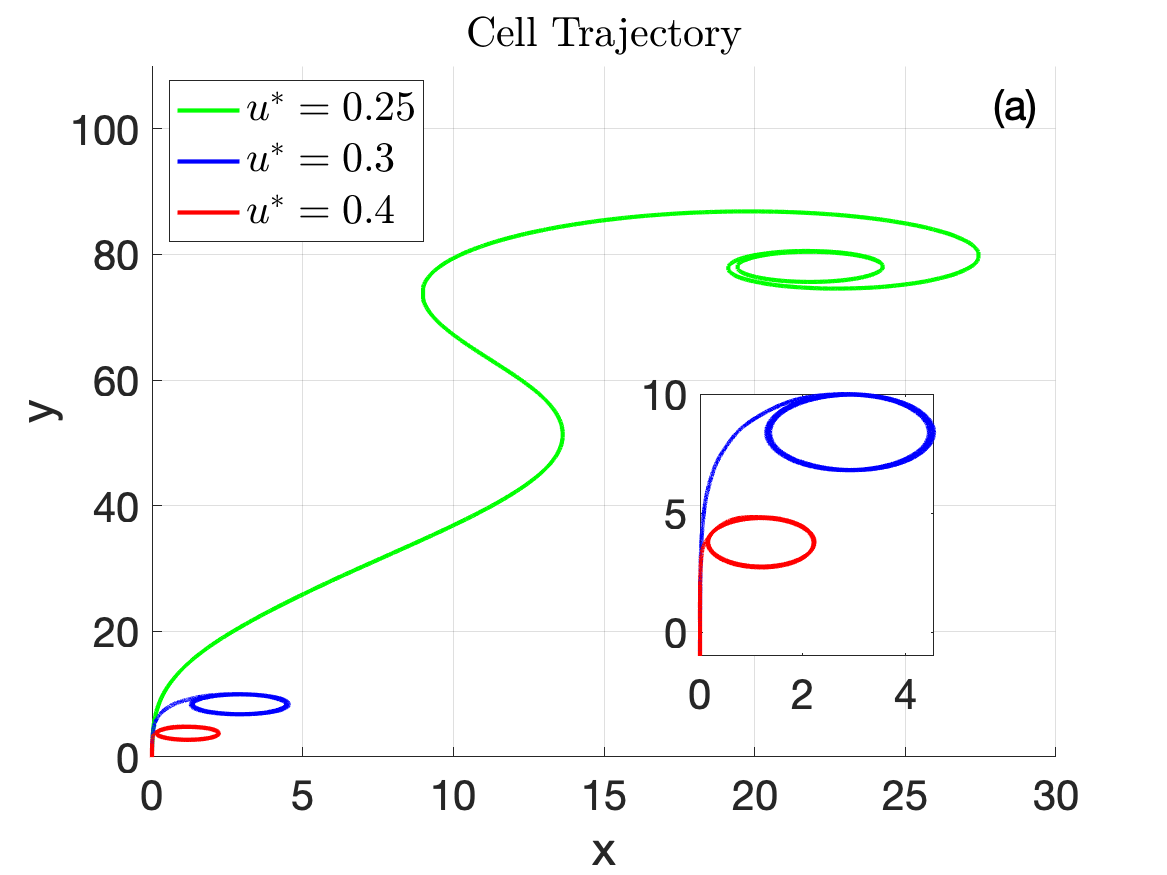}
\hspace{-20 pt}
		\includegraphics[width=0.5\linewidth]{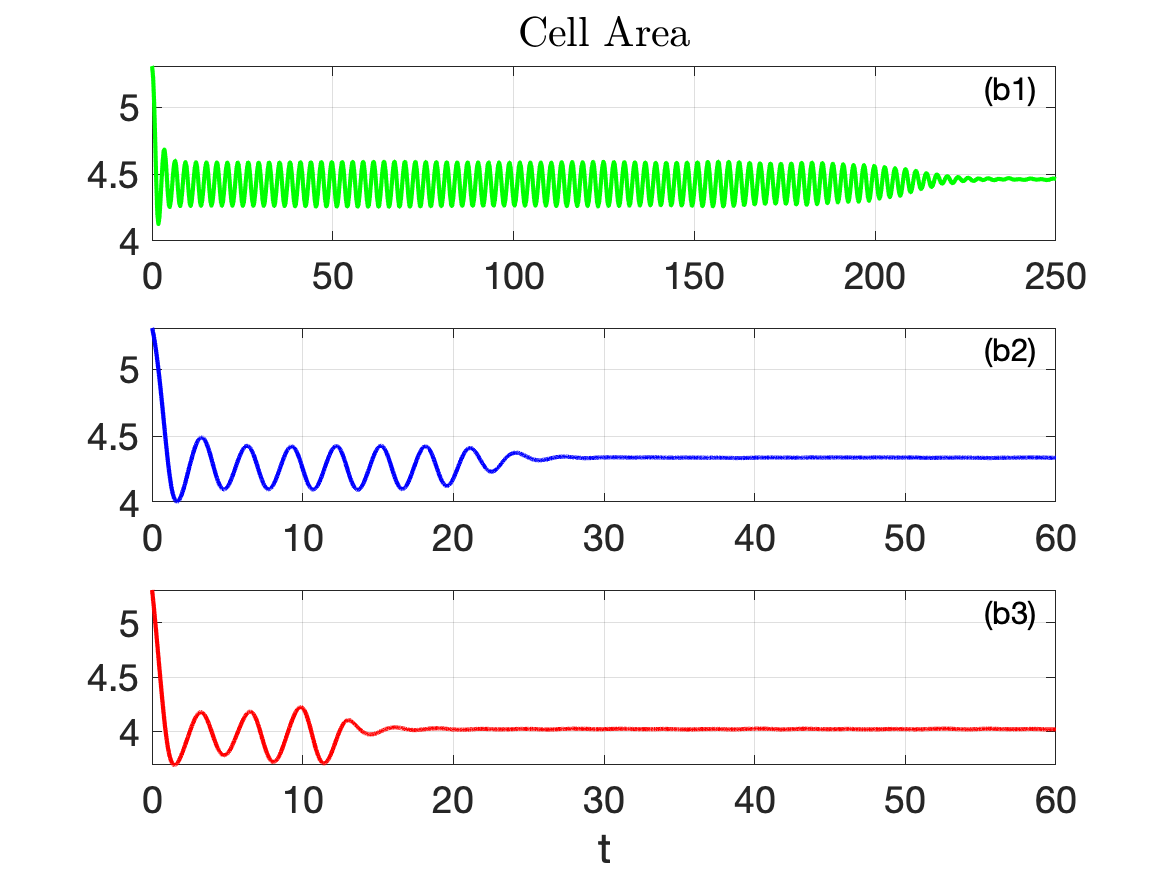}
	\caption{
(a) Cell trajectories corresponding to the three values of $u^\ast$.
(b1)--(b3) The area of the moving cell vs.\ time $t$ corresponding to the three values of $u^\ast.$
}
	\label{fig:contractivity}
\end{figure}

Figure~\ref{fig:contractivity} (a) shows the three circular trajectories corresponding 
to the three $u^\ast$ values. 
We observe that a larger value of $u^\ast$ corresponds to an earlier onset of the circular mode and 
the circle is smaller.  Figure~\ref{fig:contractivity} (b1)--(b3) show the area of the moving cell vs.\ time for 
the three sets of $u^\ast$ values as marked in Figure~\ref{fig:contractivity} (a). 
We observe again that a larger value of $u^\ast$ takes a shorter period of time before circulates. 

\medskip

\begin{figure}[H]
	\centering
		\includegraphics[width=0.45\linewidth]{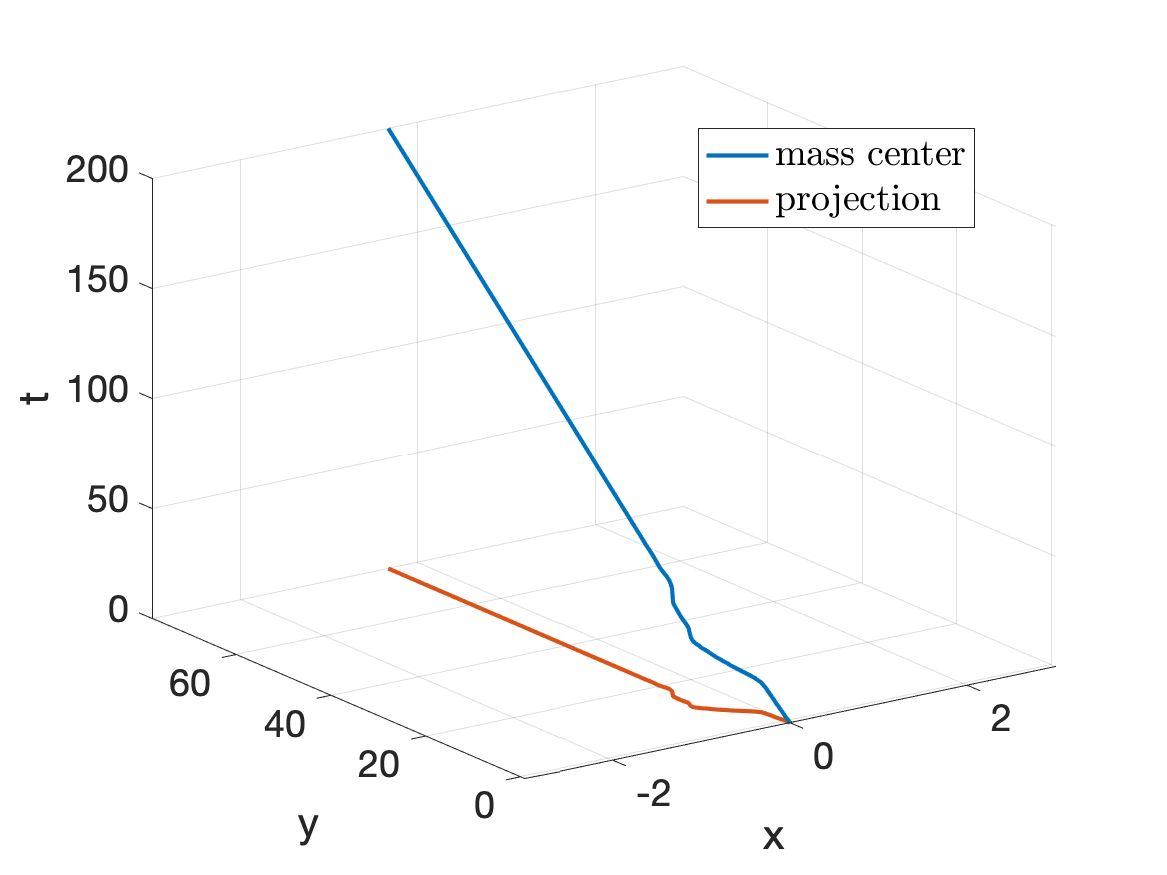}
		\includegraphics[width=0.45\linewidth]{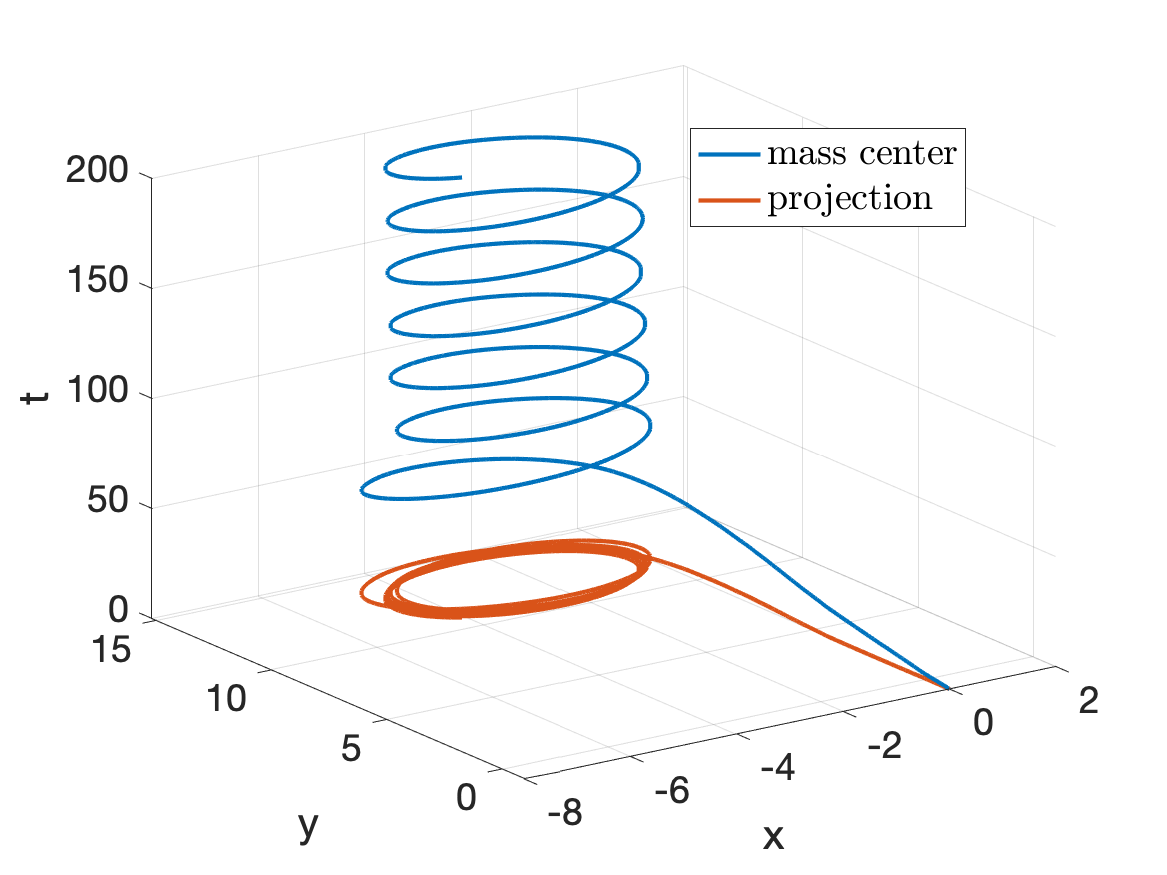}
	\caption{
With nearly identical parameters, simulations predict 
a linear trajectory with the  two-species model (left), and 
a circular trajectory with the  one-species model (right). 
}
\label{ex5:comparison}
\end{figure}

\noindent
{\bf Two-species model vs.\ one-species model.}
We simulate the cell movement with both the two-species model and the one-species model, in non-dimensionalized forms.  We use the same set of parameters including the same diffusion constant $0.1$ for the active form $u$. For the two-species model, we set the diffusion constant $D_v = 10$, this corresponds to the large-diffusion approximation in the one-species model.  While the simplified, one-species  model can can capture the leading-order behaviour of the system 
\cite{Keshet_SIAP2011}, it may be different from the full, two-species model in predicting quantitatively some of the biophysical properties of a moving cell. In Figure \ref{ex5:comparison}, we see that 
the two-species model predicts a linear trajectory (left), while the one-species model 
predicts a circular trajectory (right).

%% file: Conclusions.tex
\section{Conclusions}
\label{s:Conclusions}

We have studied the cell polarity and movement within the modelling framework of reaction-diffusion
equations and moving cell boundary. 
In particular, we have carefully examined the wave-pinning model, both the two-species
and the reduced one-species model. 

Early studies included the one-dimensional analysis of the wave pinning mechanism 
\cite{mori2008wave, Keshet_SIAP2011}
and the two-dimensional phase-filed simulation and the sharp-interface analysis
with a reduced model for cell polarization and movement 
\cite{camley2017crawling}. 
Here, we have derived the sharp-interface model
as the limit of the phase-field model as the small parameter $\ve \to 0$ with a general
two or three-dimensional setting.  Our rigorous analysis provides 
a close link between the two types of models.  

We have also developed and implemented 
a robust numerical method for the simulation of cell polarization 
and movement using the derived sharp-interface model in two-dimensional space. 
Our approach combines the level-set method for the moving cell boundary and accurate discretization 
techniques for solving the reaction-diffusion equations on the moving cell region. 
The method and algorithm pass the convergence test. 

We have done extensive numerical simulations using the full, two-species reaction-diffusion moving 
cell boundary model as well as its reduced one-species model in two space dimension. 
We find that the cell polarization is a robust process that can be triggered by various external 
stimulus with a large set of parameters, confirming the 
wave-pinning mechanism as proposed in  \cite{mori2008wave, Keshet_SIAP2011}. 
We have also traced the cell trajectory during long-time simulations. By choosing different set of 
parameters of the diffusion constants and the threshold value of the 
concentration of an active Rho GTPase protein in the normal velocity, 
we have been able to capture both the linear and circular trajectories. 
For a circular trajectory, a period of preparation time is observed. 
The full, two-species model and the reduced, single-species model predict different such 
preparation times. Therefore, the infinite diffusion of the second species which is the 
assumption of the reduce model, may need to be corrected for quantitative predictions
of different complex process of cell motility. 

In our simulations we have observed that the cell area and the $x$ and $y$ components
of the velocity at center of mass of the cell to be oscillatory 
during the period of preparation time before the cell starts to rotate completely; 
cf.\ Figure~\ref{fig:diffusionrate} and Figure~\ref{fig:contractivity}.
These are unlikely caused by numerical errors and instabilities as no such 
oscillations occur once the cell starts to rotate. 
We will investigate such oscillations further in our subsequent works. 

With our analytical tools, robust numerical methods, and computer code, we can study further the 
cell polarity and movement in several directions. 
\begin{compactenum}
\item[(1)] 
We can include many more biological components in our models and simulation. 
The first of them is the fluid flow which can be modelled by Stokes flow 
\cite{ShaoLevineRappel_PNAS12, vanderlei2011computational}. 
The boundary velocity of a cell moving around within such a flow can 
be determined by the force balance. 
The second component is the combination of attachment to and detachment from a substrate of a moving cell 
\cite{ShaoLevineRappel_PNAS12}. 
\item[(2)]
With the similar approach and simulation method, we can study the interaction and movement
of a cluster of cells, where the cell coordination and cooperation will be crucial
\cite{RappelKeshet_Rev2017,Keshet_Rev2020}. 
Likely, such studies can help understand better the molecular basis as well as mechanical forces that 
determine such an important collective biological process. 
\item[(3)]
With speeding up computations, e.g., by implementation with the GPU, 
it is possibe now for us to simulate the cell movement in a full, three-dimensional setting. 
\end{compactenum}
